\documentclass[10pt]{article}

\usepackage{verbatim,amssymb,hyperref,xcolor,amsmath,amsthm,accents,color,soul,bbm,graphicx,todonotes}
\usepackage[shortlabels]{enumitem}

\usepackage[capitalize,nameinlink,noabbrev]{cleveref}
\crefname{enumi}{item}{items}
\crefname{equation}{}{}
\crefname{subsection}{Subsection}{Subsections}
\crefname{figure}{Figure}{Figures}

\usepackage{geometry}
\hypersetup{
    colorlinks,
    linkcolor={blue!80!black},
    citecolor={green},
    urlcolor={blue!80!black}
}

\newcommand{\black}{\color[rgb]{0,0,0}}

\newcommand{\eps}{\varepsilon}

\newcommand{\dbtilde}[1]{\accentset{\approx}{#1}}

\newcommand{\cF}{\mathcal F}
\newcommand{\cC}{\mathcal C}
\newcommand{\cD}{\mathcal D}
\newcommand{\cE}{\mathcal E}

\newcommand{\cU}{\mathcal U}

\newcommand{\grad}{\nabla}

\newcommand{\sfrac}[2]{\mbox{$\frac{#1}{#2}$}}

\newcommand{\sign}{\mathrm{sign}}

\newcommand{\IB}{\mathbb B}

\newcommand{\II}{\mathbb I}

\newcommand{\1}{1\hspace{-0.098cm}\mathrm{l}}

\renewcommand{\P}{{\mathbb P}}

\newcommand{\N}{{\mathbb N}}

\newcommand{\E}{{\mathbb E}}

\newcommand{\R}{{\mathbb R}}

\newcommand{\IU}{{\mathbb U}}

\newcommand{\bas}[1]{\begin{align}\begin{split} #1 \end{split}\end{align}}

\newtheorem{theorem}{Theorem}[section]
\newtheorem{proposition}[theorem]{Proposition}
\newtheorem{lemma}[theorem]{Lemma}
\newtheorem{corollary}[theorem]{Corollary}

\newtheorem{definition}[theorem]{Definition}

\newtheorem{remark}[theorem]{Remark}

\makeatletter
\@namedef{subjclassname@2020}{%
	\textup{2020} Mathematics Subject Classification}
\makeatother

\makeatletter
\ExplSyntaxOn
\seq_new:N \g_abbrs
\prop_new:N \g_abbr_counts
\tl_new:N \l_abbr_count_tl

\ExplSyntaxOn

\bool_new:N \g_forexample

\NewDocumentCommand{\eg}{ o }{
	\IfValueT{#1}{
		\str_if_eq:noTF {fe} {#1} {
			\bool_gset_true:N \g_forexample
		} {\bool_gset_false:N \g_forexample}
	}
	\bool_if:nTF { \g_forexample } {
		\bool_gset_false:N \g_forexample
		for~example
	}{
		\bool_gset_true:N \g_forexample
		for~instance
	}
}

\NewDocumentCommand{\abbr}{m m O{#1} m m O{#4} m}{
	\expandafter\newcommand\csname#3\endcsname[1][]{
		\seq_if_in:NnTF \g_abbrs {#1} {
			\prop_get:NnN \g_abbr_counts {#1} \l_abbr_count_tl
			\prop_gput:Nnx \g_abbr_counts {#1} {\int_eval:n {\l_abbr_count_tl + 1}}
			\hyperref[#1]{#7}
		} {
			\seq_gput_left:Nn \g_abbrs {#1}
			\prop_gput:Nnn \g_abbr_counts {#1} {1}
			\expandafter\gdef\csname#1@def\endcsname{#2}
			\phantomsection\label{#1}
			\str_if_eq:nnTF{##1}{}{\emph{#2}}{##1}~(\hyperref[#1]{#7})
		}
	}
	\expandafter\newcommand\csname#6\endcsname[1][]{
		\seq_if_in:NnTF \g_abbrs {#1} {
			\prop_get:NnN \g_abbr_counts {#1} \l_abbr_count_tl
			\prop_gput:Nnx \g_abbr_counts {#1} {\int_eval:n {\l_abbr_count_tl + 1}}
			\hyperref[#1]{#4}
		} {
			\expandafter\gdef\csname#1@def\endcsname{#5}
			\seq_gput_left:Nn \g_abbrs {#1}
			\prop_gput:Nnn \g_abbr_counts {#1} {1}
			\phantomsection\label{#1}
			\str_if_eq:nnTF{##1}{}{\emph{#5}}{##1}~(\hyperref[#1]{#4})
		}
	}
}

\ExplSyntaxOff
\makeatother
\abbr{ODEs}{ordinary differential equations}{ODE}{ordinary differential equation}{ODEs}
\abbr{SDEs}{stochastic differential equations}{SDE}{stochastic differential equation}{SDEs}
\abbr{GD}{gradient descent}{GDs}{gradient descents}{GD}
\abbr{ANN}{artificial neural network}{ANNs}{artificial neural networks}{ANN}
\abbr{RMSprop}{root mean square propagation \SGD}{RMSprops}{the root mean square propagation \SGD}{RMSprop}
\abbr{Adam}{adaptive moment estimation 
}{Adams}{adaptive moment estimation \SGD}{Adam}
\abbr{Adagrad}{adaptive gradient \SGD}{Adagrads}{adaptive gradient \SGD}{Adagrad}
\abbr{SOP}{stochastic optimization problem}{SOPs}{stochastic optimization problems}{SOP}
\abbr{OP}{optimization problem}{OPs}{optimization problems}{OP}
\abbr{DNN}{deep neural network}{DNNs}{deep neural networks}{DNN}
\abbr{SGD}{stochastic gradient descent}{SGDs}{stochastic gradient descent}{SGD}
\abbr{iid}{independent and identically distributed}{i.i.d}{independent and identically distributed}{i.i.d.}
\abbr{DL}{Deep learning}{DLs}{Deep learning}{DL}
\abbr{AI}{artificial intelligence}{AIs}{artificial intelligence}{AI}
\abbr{LLM}{large language model}{LLMs}{large language models}{LLM}
\abbr{PDE}{partial differential equation}{PDEs}{partial differential equations}{PDE}
\abbr{as}{almost surely}{ass}{almost surely}{a.s.}
\abbr{pas}{$\P$-almost surely}{pass}{almost surely}{$\P$-a.s.}
\abbr{ReLUs}{rectified linear unit}{ReLU}{rectified linear unit}
\makeatother

\begin{document}

\title{ODE approximation for the Adam algorithm:\\General and overparametrized setting}

\author{Steffen Dereich$^{1}$, Arnulf Jentzen$^{2,3}$, and Sebastian Kassing$^{4}$
	\bigskip
	\\
	\small{$^1$ Institute for Mathematical Stochastics, Faculty of Mathematics and Computer Science,}\vspace{-0.1cm}\\
\small{University of M\"unster, Germany, e-mail: \texttt{steffen.dereich@uni-muenster.de}}
\smallskip
    \\
	\small{$^2$ School of Data Science and School of Artificial Intelligence, The Chinese University}
	\vspace{-0.1cm}\\
	\small{of Hong Kong, Shenzhen (CUHK-Shenzhen), China, e-mail: \texttt{ajentzen@cuhk.edu.cn}}
	\smallskip
	\\
 \small{$^3$ Applied Mathematics: Institute for Analysis and Numerics, Faculty of Mathematics and}
	\vspace{-0.1cm}\\
	\small{Computer Science, University of M{\"u}nster, Germany, e-mail: \texttt{ajentzen@uni-muenster.de}}
	\smallskip
	\\
    \small{$^4$ School of Mathematics and Natural Sciences, University of Wuppertal,}
	\vspace{-0.1cm}\\
	\small{Germany, e-mail: \texttt{kassing@uni-wuppertal.de}}
 \smallskip
	\\
}

\maketitle

\begin{abstract} 
The Adam\ optimizer is currently presumably the most popular optimization method in deep learning. 
In this article we develop an ordinary differential equation (ODE) based method to study the Adam optimizer in a fast-slow scaling regime. 
For fixed momentum parameters and vanishing step-sizes, we show that the Adam algorithm is an 
\emph{asymptotic pseudo-trajectory} of the flow of a particular vector field, 
which is referred to as the \emph{Adam vector field}. 
Leveraging properties of asymptotic pseudo-trajectories, we establish convergence results for the Adam algorithm. 
In particular, in a very general setting we show that if the Adam algorithm converges, then the limit must be a \emph{zero of the Adam vector field}, 
rather than a local minimizer or critical point of the objective function. 
	
In contrast, in the overparametrized empirical risk minimization setting, the Adam algorithm is able to locally find the set of minima. 
Specifically, we show that in a neighborhood of the global minima, the objective function serves as a Lyapunov function for the flow induced 
by the Adam vector field. As a consequence, if the Adam algorithm enters a neighborhood of the global minima infinitely often, 
it converges to the set of global minima.
\end{abstract}


\section{Introduction} 

In this article we analyze stochastic approximation algorithms for finding zeros of a vector field 
$ \mathcal{V} \colon \R^d \to \R^d $
that satisfies for all $ \theta \in \R^d $ that
\begin{equation} 
\label{eq:opt}
  \mathcal V( \theta ) = \E\bigl[ X( U, \theta ) \bigr] 
\end{equation}
where $ ( \Omega, \mathcal{F}, \P ) $ is the underlying probability space, 
where $ \mathcal{U} $ is a measurable space, 
where $ U \colon \Omega \to \mathcal{U} $ is a random variable, 
and where 
$ X \colon \mathcal U \times \R^d \to \R^d $ 
is a product measurable function which satisfies for all $ \theta \in \R^d $ that
$
  \E\bigl[ | X( U, \theta ) | \bigr] < \infty
$. 
We assume that access to $ \mathcal V $ is limited to observations 
$ ( X( U_n, \theta ) )_{ n \in \N } $ only, 
where $ ( U_n )_{ n \in \N } $ is a sequence of independent 
$ \mathcal L( U ) $-distributed random variables and where $ \theta $ can be adjusted in each step.

The first stochastic approximation algorithm for iteratively solving such problems 
was introduced in 1951 by Robbins and Monro \cite{RM51}. 
Nowadays, stochastic approximation methods have become indispensable for training \DNNs\ in machine learning tasks. 
In deep learning, the objective is typically to minimize a differentiable function $R \colon \R^d \to \R$, 
where only noisy estimates of the gradient are available, \eg, 
by evaluating the gradient on a small subset of the training data. We recover problem \eqref{eq:opt} by setting 
$ \mathcal V = - \nabla R $. 
In this context, where the vector field coincides with the negative gradient of the objective function, 
stochastic approximation algorithms for \cref{eq:opt} are usually referred to as \SGD\ optimization methods 
\cite{Ruder2016arXiv,JentzenKuckuckvonWurstemberger2025arXiv}.

Recent years have seen the development of numerous first order optimization algorithms designed to address increasingly complex and noisy problems, as well as the non-convexity of objective functions appearing in machine learning applications. Among these, algorithms such as Polyak's heavy ball \cite{polyak1964some}, Nesterov accelerated gradient \cite{nesterov1983method}, RMSProp \cite{RMSprop} and AdaGrad \cite{duchi2011adaptive} have become foundational tools for tackling large-scale optimization tasks in the presence of stochasticity. These methods typically incorporate strategies like adaptive step-sizes and momentum to accelerate convergence and improve stability. One of the most popular variants of such an adaptive optimization method is the \Adam\ optimizer proposed in \cite{kingma2015adam}, which combines the adaptive step-size approach from the RMSprop optimizer with an accelerated heavy ball method.

Despite its widespread use in practical  applications, the \Adam\ algorithm is theoretically only partially understood. 
In this work, we analyze convergence of the \Adam\ optimizer in a fast-slow scaling regime 
using the \ODE\ method; cf., \eg, \cite{ljung1977analysis, benveniste1990adaptive, benaim1996dynamical, benaim2006dynamics}.
For fixed momentum parameters and vanishing step-sizes we show that \Adam\ is an 
\emph{asymptotic pseudo-trajectory} of the semiflow \eqref{eq:ODE} for a particular vector field defined in \eqref{eq:fvector}, 
which is referred to as the \emph{\Adam\ vector field} \cite{dereich2024convergence}. 
In a second step, we use the dynamics of the semiflow to derive convergence statements for the \Adam\ optimizer depending on properties of the \Adam\ vector field. 
In particular, we show that if the \Adam\ optimizer converges, 
then the limit point is a zero of the \Adam\ vector field. 
Moreover, when the stochastic perturbation vanishes at the critical points, as is the case in the training of overparametrized \ANNs\ \cite{schmidt2013fast, vaswani2019fast, wojtowytsch2023stochastic, gess2023convergence, JiangfengLu}, we establish local almost sure convergence to the set of critical points of the original vector field $\mathcal V$.

Let us start with a rigorous definition of the \Adam\ algorithm from \cite{kingma2015adam} 
(cf., \eg, \cite[Definition~2.2]{dereich2024convergence}).

\begin{definition}[\Adam\ optimizer] 
\label{def:Adam}
Let $d\in\N$, let $ ( \cU, \IU ) $ be a measurable space, 
let $ X \colon \cU\times \R^d\to \R^d $ be product measurable, 
let $ ( \gamma_n )_{ n \in \N } $ be a sequence of strictly positive reals, 
let $ \alpha, \beta \in [0,1) $, $ \epsilon \in (0,\infty) $, $ \vartheta \in \R^d $, 
let $ ( \Omega, \mathcal{F}, \P ) $ be a probability space, 
and let $ U \colon \Omega \to \cU $ be a random variable. 

An \emph{\Adam\ optimizer (or \Adam\ algorithm) 
with innovation $ ( X, U ) $, step-sizes $ ( \gamma_n )_{ n \in \N } $,  
and damping parameters $ ( \alpha, \beta, \epsilon ) $ 
started in $ \vartheta $} 
is a stochastic process $ \theta = ( \theta_n )_{ n \in \N_0 } \colon \N_0 \times \Omega \to \R^d $ 
which satisfies that there exist
$ U_n \colon \Omega \to \cU $, $ n \in \N $, 
\iid\ copies of $ U $ such that 
for all $ n \in \N $, $ j \in \{ 1, 2, \dots, d \} $ 
it holds that $ \theta_0 = \vartheta $ and 
	\begin{equation} 
	\label{eq:Adamdefi}
	  \theta_n^{ (j) } 
	  = \theta_{ n - 1 }^{ (j) } 
	  + 
	  \gamma_n \sigma_n^{ (j) } m_n^{ (j) } ,
	\end{equation}
	where
	$m_0=v_0=0 \in \R^d$ are the starting values and where for all $ n \in \N $, $ j \in \{ 1, \dots, d \} $ 
	we have that
	\begin{align} \begin{split}\label{eq:Adamdefi2}
			m_n&= \alpha m_{n-1}+(1-\alpha) X(U_n, \theta_{n-1}), \\
			v_n^{(j)} &= \beta v_{n-1}^{(j)} +(1-\beta) (X^{(j)}(U_n, \theta_{n-1}))^2 \\
			\sigma_n^{(j)}&= \bigl( [ v_n^{(j)}/(1-\beta^n) ]^{ 1 / 2 } + \epsilon \bigr)^{-1}.
		\end{split}
	\end{align}
\end{definition}

\begin{remark}
	The definition of the \Adam\ algorithm in \cite{kingma2015adam} applies a bias correction to  the  momentum term $(m_n)_{n \in \N}$ and the second moment estimator $(v_n)_{n \in \N}$ so that the iterative description is of the form
	\begin{align} \label{eq:Adam2}
		\theta_n^{(j)} = \theta_{n-1}^{(j)} + \gamma_n (1-\alpha^n)^{-1} \sigma_n^{(j)} m_n^{(j)}.
	\end{align}
	Note that this description agrees with \cref{eq:Adamdefi,eq:Adamdefi2} 
	when employing for every $ n \in \N $ the value 
	$ \gamma_n ( 1 - \alpha^n )^{ - 1 } $
	for the step-size $ \gamma_n $ in \cref{eq:Adamdefi}.
\end{remark}

Typically, the \Adam\ algorithm is used to minimize a differentiable objective function 
$ R \colon \R^d \to \R $, where one is only able to simulate a noisy estimator 
for the gradient in the sense that 
for all $ \theta \in \R^d $ we have that
\begin{equation} 
\label{eq:gradientestimator}
  \E[ X( U, \theta ) ] = - \nabla R( \theta ) 
  .
\end{equation}
In this setting, the \Adam\ algorithm in \eqref{eq:Adamdefi} is a gradient-based optimization method 
that incorporates past gradient information through a momentum term $ ( m_n )_{ n \in \N_0 } $ 
while adaptively rescaling each coordinate using the second-moment estimator $ ( v_n )_{ n \in \N_0 } $.

In this article, we derive an \ODE\ approximation for the \Adam\ algorithm 
with fixed damping parameters $ ( \alpha, \beta, \epsilon ) $ 
and a decreasing sequence of step-sizes $ (\gamma_n )_{ n \in \N } $. 
The original work \cite{kingma2015adam} by Kingma and Ba suggests 
the default values $\alpha=0.9$, $\beta=0.999$, $\gamma=0.001$ and $\epsilon=10^{-8}$. 
Subsequent research has found that, although \Adam\ adapts its parameterwise step-sizes, 
\Adam\ fails to converge if the learning rates do not converge to zero \cite{dereichetal2024nonconvergence}
and decaying step-size schedules can substantially improve \Adam's performance \cite{wilson2017marginal, loshchilov2018decoupled}, 
which motivates analyzing \Adam\ in a setting with separating time-scales for \cref{eq:Adamdefi} and \cref{eq:Adamdefi2}.

The first approximation result in this scaling regime was given in \cite{ma2022qualitative} 
for the deterministic minimization of a differentiable function $ R \colon \R^d \to \R $ 
where one has access to the unperturbed gradient $\nabla R$, i.e., 
in the situation where one has for all $ n \in \N $ that 
$ X( U_n, \theta_{ n - 1 } ) = ( \nabla R )( \theta_{ n - 1 } ) $. 
Therein, it is shown that the \Adam\ optimizer 
can be approximated by the solution to the \ODE 
\begin{align} \label{eq:Eflow}
  \dot \vartheta_t = -\frac{ ( \nabla R )( \vartheta_t ) }{ | ( \nabla R )( \vartheta_t ) | + \epsilon } .
\end{align}
Since $ \epsilon $ is a small hyperparameter, 
the dynamics of \eqref{eq:Eflow} are similar to 
signGD $ \dot \vartheta_t = -\sign( ( \nabla R )( \vartheta_t ) ) $. 
Moreover, the solution of \eqref{eq:Eflow} satisfies 
$ \frac{ d }{ dt } R( \vartheta_t ) = - | ( \nabla R )( \vartheta_t ) |^2 / ( | ( \nabla R )( \vartheta_t ) + \epsilon | ) $ 
so that the \ODE\ approximation typically implies 
convergence of the \Adam\ optimizer to the set of critical points of $ R $.

In the stochastic setting, the situation is much more intricate. Here, we draw our intuition 
from the analysis of random fast-slow systems; see, \eg, \cite[Chapter 7]{freidlin2012random}. 
Stochastic approximation schemes with two separating time-scales have already 
been studied in the scientific literature, most prominently by Borkar~\cite{borkar1997stochastic, borkar2008stochastic,  borkar2025stochastic}. 
A typical assumption in these references is that the step-sizes converge to zero and, at the same time, the momentum parameters converge to one 
but on different time-scales, i.e. $1-\alpha_n \sim 1- \beta_n$ and  $\gamma_n = o (1-\alpha_n)$.
This asymptotic regime is fundamentally different from the regime where the momentum parameters $\alpha$ and $\beta$ converge to $1$ with the same speed as the step-size converges to zero, which leads to a second order limiting equation similar to Polyak's heavy ball; 
see \cite{barakat2021convergence, malladi2022sdes}.

In this work, we follow the line of research on stochastic approximation schemes with separating time-scales and derive a homogenization result. 
Since the momentum parameters for the \Adam\ algorithm are typically fixed hyperparameters that do not change over time, we consider the small step-size limit for fixed $(\alpha, \beta, \epsilon)$. Hence, the iterates defined by \eqref{eq:Adam2} preserve their discrete character asymptotically.

In our main statement in \cref{thm:main} below we show that the \Adam\ algorithm 
can be approximated by the semiflow of functions $ \Psi_t \colon \R^d \to \R^d $, $ t \in [0,\infty) $, 
solving the \ODE\ that for all 
$ \theta \in \R^d $, 
$ t \in [0,\infty) $
it holds that
\begin{equation} \label{eq:ODE}
  \dot \Psi_t( \theta ) = f( \Psi_t( \theta ) )
  \qquad 
  \text{and}
  \qquad 
  \Psi_0(\theta)=\theta
\end{equation}
where 
$ f \colon \R^d \to \R^d $ 
is the so called \Adam\ vector field that has been introduced in \cite{dereich2024convergence} and that we recall here in this work within the next notion.

\begin{definition}[\Adam\ vector field]
\label{def:Adam_vector_field}
Let $ \alpha, \beta \in [0,1) $, $ \epsilon \in (0,\infty) $ 
and let $ ( X, U ) $ 
be an innovation in the sense of Definition~\ref{def:Adam}. 
Then we say that $ f $ is the 
\emph{\Adam\ vector field of the innovation $ ( X, U ) $ 
with damping parameters $ ( \alpha, \beta, \epsilon ) $}
if and only if it holds that 
$ f = ( f^{ (1) }, \dots, f^{ (d) } ) \colon \R^d \to \R^d $ 
is the function which satisfies 
for all $ \theta \in \R^d $, $ j \in \{ 1, 2, \dots, d \} $ that
\begin{equation} 
\label{eq:fvector}
  f^{ (j) }( \theta ) 
  = 
  ( 1 - \alpha ) 
  \,
  \E\Biggl[
    \biggl[
      \Bigl[ 
        ( 1 - \beta )
\textstyle 
        \sum_{ k = -\infty }^0 
        \beta^{ - k } 
        | X^{ ( j ) }( U_k, \theta ) |^2
	  \Bigr]^{ 1 / 2 } 
	  + \epsilon 
	\biggr]^{ - 1 } 
\displaystyle 
	\sum_{ k = - \infty }^0 
	\alpha^{ - k } 
	X^{ (j) }( U_k, \theta ) 
  \Biggr] 
  .
\end{equation}
\end{definition}

To prove the \ODE\ approximation we work with the following regularity assumption 
for the innovation $ ( X, U ) $.

\begin{definition}[Locally regular innovation]  
\label{assu:regular}
Let $ p \in [2,\infty) $. Then we say that the innovation $ ( X, U ) $ 
is \emph{locally $ p $-regular} if and only if there exists $ c \in C( \R, \R ) $ such that 
for all $ \theta, \vartheta \in \R^d $ it holds that
\begin{equation} 
\label{eq:34296255544}
  \E\bigl[ | X(\theta,U) |^p \bigr]
  \le  c(0) ( 1 + e^{ c(0) |\theta| } )
\qquad 
  \text{and}
\qquad 
  \E\bigl[ 
    | X( U, \theta ) - X( U, \vartheta ) |^p 
  \bigr]
  \le c( |\theta| + | \vartheta | ) | \theta - \vartheta |^p 
  .
\end{equation}
\end{definition}
Throughout this work, we assume that $ \alpha < \sqrt{\beta} $, which is, \eg, satisfied 
in the case of the default values $ \alpha = 0.9 $ and $ \beta = 0.999 $. 
Under this condition, local $ p $-regularity of the innovation $ ( X, U ) $ implies 
global boundedness and local Lipschitz continuity of the \Adam\ vector field $ f $; 
see \cite[Lemma~3.1 and (148)]{dereich2024convergence}. Therefore, 
the \ODE\ defining the semiflow $ ( \Psi_t )_{ t \in [0,\infty) } $ in \cref{eq:ODE} 
admits a unique solution. 

In our main result in \cref{thm:main} below we establish that, asymptotically, 
the dynamics of the \Adam\ algorithm are governed 
by the flow of the \Adam\ vector field in \eqref{eq:fvector}. 
This stands in sharp contrast to \SGD, whose limiting dynamics are driven by the mean vector field 
$ \mathcal V $ 
in \cref{eq:opt} 
(see, \eg, \cite{benaim2006dynamics, mertikopoulos2020almost}). 
In general, the \Adam\ vector field has different zeros and might point in a different direction compared to the mean vector field $ \mathcal V $. 
Moreover, even if the mean vector field is given as the gradient of a differentiable function, as in \eqref{eq:gradientestimator}, 
the \Adam\ vector field \cref{eq:fvector} need not possess a gradient-like structure. 
Consequently, the \Adam\ algorithm may fail to converge.

To state the main theorem that relates the \Adam\ algorithm to the semiflow of the \Adam\ vector field, 
we introduce an appropriate time scaling. 
Specifically, given a sequence of step-sizes $ ( \gamma_n )_{ n \in \N } \subseteq (0,\infty) $ 
we define the associated \emph{training times} 
$ ( t_n )_{ n \in \N_0 } \subseteq \R $ 
by $ t_n = \sum_{i=1}^n \gamma_i$ for all $ n \in \N_0 $.

\begin{theorem}[\ODE\ approximation for \Adam] 
\label{thm:main}
Let $ \epsilon \in (0,\infty) $, $ \alpha \in [0,1) $, $ \beta \in ( \alpha^2, 1 ) $, 
let $ ( \gamma_n )_{ n \in \N } \subseteq (0,\infty) $ be a non-summable decreasing zero sequence satisfying
	\bas{\label{eq:regular_steps}
	\textstyle  
		\sup_{ \varepsilon \in (0,\infty) }
		\liminf_{ n \to \infty } \bigl( ( \gamma_n )^{ - 1 } \gamma_{ n + \lceil ( \gamma_n )^{ - 1 } \varepsilon \rceil} \bigr) > 0 ,
	}
	let $ ( t_n )_{ n \in \N_0 } $ be the respective training times, let $ p \in [2,\infty) $, $ d \in \N $, $ \vartheta \in \R^d $, 
	let $ ( X, U ) $ be a locally $ p $-regular innovation, and 
	let $ ( \theta_n )_{ n \in \N_0 } $ be the \Adam\ algorithm with innovation $ ( X, U ) $, 
	step-sizes $ ( \gamma_n )_{ n \in \N_0 } $,  
	and 
	damping parameters $ ( \alpha, \beta, \epsilon ) $
	started in $ \vartheta $.
	Then
	\begin{enumerate}[(i)]
		\item  
		\label{thm:item_i}
		it holds for all $ T, R \in (0,\infty) $ that 
		\begin{equation*}
		  \textstyle 
		  \sup_{ 
		    n \in \N_0, \, 0 \le t_n - t_N \leq T 
		  }
		  \bigl(
		    \left| 
		      \theta_n - \Psi_{ t_n - t_N }( \theta_N ) 
		    \right| 
  		    \1_{ \{ | \theta_N | \le R \} } 
		  \bigr)
		  \overset{ N \to \infty }{ \longrightarrow } 0, 
		  \quad \text{in probability,}
		\end{equation*}
		and
		\item
		\label{thm:item_ii}
		 it holds for all $ T, R \in (0,\infty) $  with $ \sum_{ n  = 1 }^{ \infty } ( \gamma_n )^{ 1 + p/3 } < \infty $ that 
		\begin{equation*}
		  \textstyle 
		  \sup_{ 
		    n \in \N_0, \, 0 \le t_n - t_N \leq T 
		  }
		  \bigl(
		    \left| 
		      \theta_n - \Psi_{ t_n - t_N }( \theta_N ) 
		    \right| 
  		    \1_{ \{ | \theta_N | \le R \} } 
		  \bigr)
		  \overset{ N \to \infty }{ \longrightarrow } 0, 
		  \quad \text{almost surely.}
		\end{equation*}
	\end{enumerate}
\end{theorem}

We prove \cref{thm:main} in \cref{sec:proof} below. 
From a somewhat different perspective, we also refer, for example, to \cite{MR3948080,LiTaiE2015arXiv,MR4795573,gess2024stochastic,gess2024stochastic2,Callisti2025arXiv} (see also \cite{MR4253710}) for works that, on finite time intervals, investigate probabilistically weak time-continuous approximations of time-discrete \SGD\ optimization methods (in some cases including the \Adam\ optimizer \cite{Callisti2025arXiv}) via suitable \emph{stochastic modified equations} \cite{MR2783188,MR2214851}.
In the aforementioned references, the time-continuous approximation processes are solutions to \SDEs\ with the \emph{negative gradient of the objective function} 
appearing in the drift coefficients of the \SDEs. 
In contrast, \cref{thm:main} establishes a purely deterministic \ODE\ as the time-continuous limit, where the drift coefficient is not the negative gradient of the objective function but instead the \emph{\Adam\ vector field}.

\begin{remark}
In \eqref{eq:regular_steps}, we impose an assumption that controls how much the step-sizes vary over a given comparison horizon. 
Clearly, \eqref{eq:regular_steps} is satisfied if 
we choose $ ( \gamma_n )_{ n \in \N } = ( c n^{ - \gamma } )_{ n \in \N} $ 
with $ c > 0 $ and $\gamma \in (0,1] $, 
since
\begin{align*}
  \gamma_{ n + \lceil\gamma_n^{-1}\rceil} 
  \ge 
  \gamma_{ n + 2 c^{-1} n^{\gamma} } 
  \ge  
  c ( 1 + 2 c^{-1} )^{ - \gamma } 
  n^{ - \gamma } 
  = ( 1 + 2 c^{ - 1 } )^{ - \gamma } 
  \gamma_n .  
\end{align*}
We add that if the innovation $ ( X, U ) $ is not only $ p $-regular but satisfies the stronger assumption 
that there exists $ c \in \R $ such that for all $ \theta, \vartheta \in \R^d $ it holds that
$ 
  \E\bigl[ | X( U, \theta ) |^p \bigr]^{ 1 / p } \le c 
$ 
and 
$ 
  \E\bigl[ | X( U, \theta) - X( U, \vartheta )|^p \bigr]^{ 1 / p } \le c | \theta - \vartheta |
$, 
then the conclusion of \cref{thm:main} also holds 
without the indicator function $ \1_{ \{ | \theta_N | \le R \} } $ 
in \cref{thm:item_i,thm:item_ii}. 
\end{remark}

\begin{remark}\label{rem:2settings}
In what follows, we apply \cref{thm:main} 
to investigate the asymptotic behavior of the \Adam\ algorithm. 
To this end, we distinguish between two frameworks. 
\begin{itemize}
\item The setup of \cref{thm:item_i} in \cref{thm:main} 
will be referred to as \emph{stochastic convergence setup}
or setup of \cref{thm:item_i}. 
\item 
The setup of \cref{thm:item_ii} in \cref{thm:main}, 
in which we additionally assume that 
$ \sum_{ n = 1 }^{ \infty } ( \gamma_n )^{ 1 + p / 3 } < \infty $, 
will be referred to as \emph{almost sure convergence setup} 
or setup of \cref{thm:item_ii}. 
\end{itemize}
Throughout this article, $ f $ will always denote the \Adam\ vector field defined in \eqref{eq:fvector}.
\end{remark}

In the dynamical systems literature, the ODE method is employed to establish convergence results for stochastic algorithms by exploiting properties of the limiting semiflow. This approach relies on the notion of asymptotic pseudo-trajectories and the fact that, for precompact asymptotic pseudo-trajectories, the limit set is internally chain transitive; see \cite[Theorem~5.7]{benaim2006dynamics}.

\begin{definition}
	Let $(y_t)_{t \in [0,\infty)}$ be a curve in $\R^d$. We say that $(y_t)_{t \in [0,\infty)}$ is an \emph{asymptotic pseudo-trajectory} of the semiflow \eqref{eq:ODE} if and only if it holds 
	for all $ T \in (0,\infty) $ that
	$$
	\lim_{t \to \infty} \sup_{0\le h \le T} |y(t+h)-\Psi_h(y(t))|=0.
	$$ 
\end{definition}

Using the training time sequence $(t_n)_{n \in \N_0}$, we define 
a stochastic process $ \bar \theta \colon [0,\infty) \times \Omega \to \R^d $ 
by letting for every $n\in\N_0$,  $\bar \theta_{t_n}=\theta_{n} $ 
and, on each interval $[t_{n-1},t_{n}]$ with $n\in\N$,  letting  $(\bar \theta_t)_{t \in [t_{n-1},t_{n}]}$ be 
the linear interpolation between the two boundary values.
 
\cref{thm:main}, together with the boundedness of the \Adam\ vector field, 
yields the following corollary. 
Since the connection between asymptotic pseudo-trajectories 
and chain transitive sets requires precompactness of the paths, 
in the corollary we restrict attention to the event
$$
  \mathbb B = \Big\{\limsup_{n \to \infty} |\theta_n|<\infty\Big\}.
$$

\begin{corollary} \label{cor:APT}
Consider the setup of \cref{thm:item_ii} in \cref{thm:main}. 
Then,  on $\mathbb B$, the interpolated process $(\bar \theta_t)_{t \in [0,\infty)}$ as introduced above is almost surely an asymptotic pseudo-trajectory of the semiflow~\eqref{eq:ODE}.
\end{corollary}

\begin{proof}
	The proof follows directly from \cref{thm:main} and two observations. 
	First, note that  
	\begin{align} \label{eq:localarg}
		\mathbb B = \bigcup_{R \in \N} \bigcap_{n\in \N_0} \{ |\theta_n| \le R\},
	\end{align}
	so it suffices to establish the claim on the event  
	\[
	\bigcap_{n\in \N_0}\{ |\theta_n| \le R \},
	\]
	for an arbitrarily fixed $R\in\N$. 
	Second, there exists a real number $ c \in (0,\infty) $ such that for all $n \in \N$ it holds that 
	\[
	|\theta_n-\theta_{n-1}| \le c \gamma_n
	\]
	(see \cite[Lemma~3.1]{dereich2024convergence}). 
	In particular, $|\theta_n-\theta_{n-1}|\to 0$ almost surely. 
	Combining this fact and the boundedness and the Lipschitz continuity of $ f $ 
	(see (148) in \cite{dereich2024convergence}) 
	with \cref{thm:main} completes the proof of 
	Corollary~\ref{cor:APT}. 
\end{proof}

Let us mention a few implications of \cref{thm:main} for the asymptotic behavior of the \Adam\ algorithm, 
depending on the properties of the \Adam\ vector field. 
Our first corollary states that, if the \Adam\ algorithm converges, 
then its limit point has to be a zero of the \Adam\ vector field. 

\begin{corollary}[Convergence to zeros of the \Adam\ vector field]
\label{cor:criticalpoint}
Consider the setup of \cref{thm:item_i} of \cref{thm:main}, 
let $ \theta_\infty $ be an $ \R^d $-valued random variable 
and let $\Omega_0$ be an event with $\lim_{n\to\infty} \theta_n=\theta_\infty$, in probability, on $\Omega_0$, i.e., 
assume for all $ \eps \in (0,\infty) $ that 
$$
\textstyle 
  \lim_{n\to\infty} \P\bigl(\{|\theta_n-\theta_\infty|\ge \eps\}\cap \Omega_0\bigr)=0.
$$
Then 
$$
  \P\bigl(\{f(\theta_\infty)\not=0\}\cap\Omega_0\bigr) = 0 .
$$
\end{corollary}

\begin{proof}
	For $R\in(0,\infty)$ let $\Omega_0^R=\Omega_0\cap \{|\theta_\infty|<R/2\}$ and consider for $\delta\in(0,\infty)$ the mapping
	$$
		n_\delta:\N_0\to \N_0,\ N\mapsto n_\delta(N)=\inf \{n\in\N_0: t_n-t_N\ge \delta\}.
	$$
	Note that $\lim_{N\to\infty} \P(\{|\theta_N|>R\}\cap \Omega_0^R)=0$ so that by \cref{thm:main} we have that 
	\bas{\label{eq72385}
		\lim_{N\to\infty} \1_{\Omega_0^R} (\theta_{n_{\delta}(N)}-\Psi_{t_{n_\delta(N)}-t_N}(\theta_N))=0,\text{ \ in probability}.
	}
	By the uniform boundedness of $f$ and the fact that $t_{n_\delta(N)}-t_N\to \delta$, we get that
	$$
	\lim_{N\to\infty} \Psi_{t_{n_\delta(N)}-t_N}(\theta_N)-\Psi_{\delta}(\theta_N)=0.
	$$
	Recalling that $\lim_{N\to\infty} \1_{\Omega_0^R} (\theta_N-\theta_\infty)=0$, in probability and that $\Psi_\delta$ is continuous we conclude that
	$$
	\lim_{N\to\infty} \1_{\Omega_0^R} (\theta_{n_{\delta}(N)}-\Psi_{t_{n_\delta(N)}-t_N}(\theta_N)) = \1_{\Omega_0^R} (\theta_\infty -  \Psi_{\delta}(\theta_\infty)), \text{ \ in probability}.
	$$
	Consequently, we get with \cref{eq72385} and the fact that $\Omega_0=\bigcup_{R\in\N}\Omega_0^R$ that
	$$
	\P(\{\Psi_\delta(\theta_\infty)\not= \theta_\infty\}\cap \Omega_0)=0.
	$$
	This statement is true for every $\delta\in(0,\infty)$. Since $\mathbb Q\cap(0,\infty)$ is countable, the statement is also true when taking the union of all these events over $\delta\in \mathbb Q\cap(0,\infty)$ and by continuity of $\Psi_{ ( \cdot ) }(\theta)$ for every fixed $\theta$ we get that actually the statement is true when taking the union over all $\delta\in(0,\infty)$, i.e.
	$$
	\P\Bigl(
	  \bigl( 
  	  \cup_{\delta>0} \{\Psi_\delta(\theta_\infty)\not= \theta_\infty\}
  	  \bigr) \cap \Omega_0\Bigr)=0.
	$$
	This finishes the proof since this shows that the flow is almost surely constant when started in~$\theta_\infty$ on the event $\Omega_0$.
\end{proof}

Corollary~\ref{cor:criticalpoint} can be viewed as a convergence result for the \Adam\ optimizer. 
In the scientific literature several error analyses and convergence results for \Adam\ and closely related optimization methods 
have already been established, \eg, in \cite{dereich2024convergence,ding2025adam_decoupled_wd_tmlr,LiRakhlinJadbabaie2023arXiv,ReddiKaleKumar2019arXiv,XHLK24,YushunZhangetal2022arXiv} 
and the references therein. 
These results are, however, of a 
fundamental different nature than Corollary~\ref{cor:criticalpoint} above. 
Specifically, in the above named references several additional assumptions on the objective function and the gradient/innovation 
-- such as appropriate $ L $-smoothness 
(cf., \eg, \cite[Definition~1]{ZhangHeSraJadbabaie2023arXiv}, \cite[Definition~3.2]{LiRakhlinJadbabaie2023arXiv}, and \cite[Section~3.1]{LiQianTianRakhlinJadbabaie2023arXiv}) 
and suitable growth bounds on the noise (cf., \eg, \cite[Assumption~2.2]{YushunZhangetal2022arXiv}) -- are imposed 
(much beyond the hypotheses in the setup of \cref{thm:item_i} in \cref{thm:main} 
but then an upper bound for the regret/error of \Adam\ is provided and/or convergence of \Adam\ is concluded \cite{dereich2024convergence}. 
Corollary~\ref{cor:criticalpoint} is of a completely different nature as it only imposes very mild and general conditions 
on the objective function and the gradient/innovation (cf.\ \cref{assu:regular} above) 
but instead \emph{assumes} that one knows already that \Adam\ converges 
in probability to a random variable (at least on some event of the underlying probability space) but concludes then that 
this limiting random variable must be \emph{a zero of the \Adam\ vector field} in (8) above, 
rather than a zero of the gradient of the objective function. 


Next, we establish convergence to an asymptotically stable point of the \Adam\ vector field, 
assuming that  $(\theta_n)_{n \in \N_0}$ enters its basin of attraction infinitely often. 
This result extends the convergence statement in~\cite{dereich2024convergence} to a larger class of vector fields and step-sizes, 
but it does not provide a convergence rate.

\begin{corollary}
	Consider the setup of \cref{thm:item_ii} in \cref{thm:main},
	let $\mathcal D \subset \R^d$ be compact 
	and let $\theta^* \in \operatorname{int}(\mathcal D)$ satisfy 
	$$
	\lim_{t \to \infty} \Psi_t(\theta)=\theta^*, \quad \text{ uniformly in } \theta \in \mathcal D.
	$$
	Then, one has almost surely on the event
	$
	\{\theta_n \in \mathcal D \text{ for infinitely many } n \in \N_0 \} 
	$
	that
	$$
	\lim_{n \to \infty} \theta_n =\theta^*.
	$$
\end{corollary}

\begin{proof}
	Let $\eps\in(0,\infty)$ be such that $\overline {B_{2\eps}(\theta^*)}\subset \cD$. The statement follows once we showed that almost surely on the event
	$ \{\theta_n \in \mathcal D \text{ for infinitely many } n \in \N_0 \}$ one has that $(\theta_n)_{n\in\N_0}$ takes values in $\overline{B_{2\eps}(\theta^*)}$ for all but finitely many $n\in\N_0$.
	
	To see this  choose $T\in[0,\infty]$ such that for all $\theta\in\cD$ and $t\in[T,\infty)$ one has $\Psi_t(\theta)\in B_\eps(\theta^*)$. Since $\cD$ is compact we have that, almost  surely,
	$$
	\1_{\{\theta_n\in\cD\}} \sup_{m\in\N_0:0\le t_m-t_n\le 3T} |\theta_m-\Psi_{t_m-t_n}(\theta_n)|\le \eps
	$$
	is satisfied for all but finitely many $n\in\N_0$, say for $n\ge n_0(\omega)$. Choose $n_1\in\N$ such that $\gamma_{n_1}\le T$. Next, we will show by contradiction that whenever $\theta_n(\omega)\in\cD$ for a $n\in\N_0$ with $n\ge n_0(\omega)\vee n_1$ one has that for all $m\in\N$ with $t_m-t_n\ge T$, $\theta_m\in \overline{B_{2\eps}(\theta^*)}$. Suppose this were not true and denote by  $m\in\N$ the minimal index with $t_m-t_n\ge T$ and  $\theta_m(\omega)\not \in\overline { B_{2\eps}(\theta^*)}$. 
	We distinguish two cases: if $t_{m}-t_n\le 3T$ then we get a contradiction via
	$$
	|\theta_m-\theta^*|\le |\theta_m-\Psi_{t_m-t_{n}}(\theta_{n}(\omega))| + |\Psi_{t_m-t_{n}}(\theta_{m'}(\omega))-\theta^*|\le \eps +\eps.
	$$
	If $t_{m}-t_n> 3T$, then  we pick $m'\in \N$ with $t_{m'}-t_n\ge T$ and $t_m-t_{m'}\in [T,3T]$. By choice of $m$, we have that  $\theta_{m'}\in \overline{B_{2\eps}(\theta^*)}$ which yields a contradiction via
	$$
	|\theta_m-\theta^*|\le |\theta_m-\Psi_{t_m-t_{m'}}(\theta_{m'})| + |\Psi_{t_m-t_{m'}}(\theta_{m'})-\theta^*|\le \eps +\eps.
	$$
\end{proof}

The following corollary establishes a general convergence result for the semiflow \cref{eq:ODE}, 
assuming the existence of a Lyapunov function. It is obtained through 
a slight modification of \cite[Proposition~6.4]{benaim2006dynamics}.

\begin{corollary} \label{cor:omegalimit}
Consider the setup of \cref{thm:item_ii} in \cref{thm:main}, 
let $\Lambda \subset \mathcal D \subset \R^d$ be closed sets that are invariant under the semiflow~\eqref{eq:ODE}, i.e., 
for all $ t \in [0,\infty) $ one has $\Psi_t(\Lambda) \subset \Lambda$ and $\Psi_t(\mathcal D) \subset \mathcal D $. 
Let $\varphi \colon \mathcal D \to \R$ be a Lyapunov function for $\Lambda$, i.e.,
a continuous function such that $[0,\infty)\ni t \mapsto \varphi(\Psi_t(\theta))$ is constant for every $\theta \in \Lambda$ 
and strictly decreasing for every $\theta \in \mathcal D \setminus \Lambda$. Assume that for every compact $K \subset \R^d$ the set $\varphi(\Lambda \cap K)$ has empty interior. 
Then, one has for almost all $$
	\omega \in \IB \cap \{\theta_n \in \mathcal D \text{ for all but finitely many } n\}
	$$
	that the set of accumulation points
	$$
	\mathcal L\bigl((\theta_n(\omega))_{n \in \N_0}\bigr) := \bigcap_{n \in \N_0} \overline{(\theta_k(\omega))_{k \ge n}}
	$$
	is a connected subset of $\Lambda$ and $\varphi$ is constant on $\mathcal L((\theta_n(\omega))_{n \in \N_0})$.
\end{corollary}

For example, $\R^d$ and $f^{-1}(\{0\})$ are invariant for the semiflow \eqref{eq:ODE}. 
Typically, one chooses $\Lambda$ to be the set of equilibria, i.e., $\Lambda = f^{-1}(\{0\})$.

\begin{proof} By Corollary~\ref{cor:APT}, on $\IB$, 
the interpolated version $(\bar \theta_t)_{t\in[0,\infty)}$ is an asymptotic pseudo-trajectory of the semiflow.  
We fix an $\omega$ for which this is true and also fix 
an $ R \in (0,\infty) $ such that 
$(\theta_n(\omega))_{n \in \N_0} \subset B_R(0)$. 
Since $\mathcal L((\theta_n(\omega))_{n \in \N_0})$ is 
invariant for the semiflow \cref{eq:ODE} 
(see \cite[Theorem~0.1]{benaim1996asymptotic}), 
one has 
$|\Phi_t(\theta)|\le R$ for all $ t \ge 0 $ and 
accumulation point $\theta \in \mathcal L((\theta_n(\omega))_{n \in \N_0})$. Therefore, for every accumulation point the mapping $t \mapsto \varphi(\Phi_t(\theta))$ is non-increasing and bounded. The proof now follows as in~\cite[Proposition~6.4]{benaim2006dynamics}.
\end{proof}

Let us consider the situation where \Adam\ is applied to minimize a differentiable objective function 
$ R \colon \R^d \to \R $, 
under the assumptions that \cref{eq:gradientestimator} holds. 
Our aim is to apply Corollary~\ref{cor:omegalimit} to establish convergence of the \Adam\ algorithm to the set 
of critical points of $ R $. For this we interpret $ R $ as a Lyapunov function 
for the critical points $ ( \nabla R )^{ - 1 }( \{ 0 \} ) $. 
Achieving this, however, requires a careful comparison between the \Adam\ vector field \eqref{eq:fvector} 
and the negative gradient vector field $-\nabla R$. We first present the general statement, 
and in the next section we prove a refined version in the context of empirical risk minimization 
with overparametrized neural networks.

\begin{corollary} 
\label{cor:Lyapunov}
Consider the setup of \cref{thm:item_ii} in \cref{thm:main}, 
let $ R \in C^d( \R^d, \R ) $, let $\mathcal D$ be a closed set that is invariant for the semiflow \cref{eq:ODE} 
and assume for all $\theta \in \mathcal D$ that 
	\begin{align}\label{eq:Lyapunov}
		\langle f(\theta), \nabla R(\theta) \rangle < 0 \quad \text{ or } \quad  f(\theta)= \nabla R(\theta)=0.
	\end{align}
	Then one has for almost all
	$$
	\omega \in \IB \cap \{\theta_n \in \mathcal D \text{ for all but finitely many } n\}
	$$
	that the set of accumulation points
	$$
	\mathcal L\bigl((\theta_n(\omega))_{n \in \N_0}\bigr) := \bigcap_{n \in \N_0} \overline{\{\theta_k(\omega): k \ge n\}}
	$$
	is a connected subset of the critical points of $R$ and $R(\mathcal L((\theta_n(\omega))_{n \in \N_0}))$ is constant.
\end{corollary}

\begin{proof}
By Sard's theorem, the set of critical values of $ R $ is closed and has Lebesgue measure $ 0 $. 
Moreover, the objective function $ R $ is clearly a Lyapunov function for the set $ \mathcal D \cap (\nabla R)^{-1}( \{ 0 \} ) $ 
under the semiflow \eqref{eq:ODE}. The result now follows from Corollary~\ref{cor:omegalimit}. 
\end{proof}

\begin{remark}
Suppose that there exists an $ \ell \in \R $ 
such that for all $ \theta \in \mathcal D_\ell := \{\theta' \in \R^d \colon R(\theta') \le \ell\} $ 
statement \cref{eq:Lyapunov} is true. 
Then $ \mathcal D_\ell $ is an invariant set for the semiflow \cref{eq:ODE}. 
Indeed, for every fixed $\theta\in \mathcal D_\ell$ one has
	\bas{ \frac d{dt} R(\Psi_t(\theta))=\langle f(\Psi_t(\theta)), \nabla R(\Psi_t(\theta))\rangle \le 0
	}
	on the time interval $\II:=[0,\tau]\cap[0,\infty)$, where the exit time $\tau := \inf\{t \ge 0: \Psi_t(\theta) \notin \mathcal D_\ell\}$. Hence, we have by closedness of $\cD_\ell$ that on $\II$, $(\Psi_t(\theta))$ stays in $\cD_\ell$.
	If $\tau$ were finite, then either $f(\Psi_\tau(\theta))=0$ or $f(\Psi_\tau(\theta))\not=0$. In the first case this would imply that $\Psi_t(\theta)=\Psi_\tau(\theta)$ for all $t\ge \tau$ and in the second case  one has 
		\bas{\frac d{dt} R(\Psi_\tau(\theta))=\langle f(\Psi_\tau(\theta)), \nabla R(\Psi_\tau(\theta))\rangle < 0
	}	
	so that in both cases  $\tau$ cannot be an exit time.
\end{remark}

\section{Degenerate Noise and Empirical Risk Minimization}
In this section, we consider the \Adam\ algorithm in the context of stochastic optimization that is the case where $\E[X(U,\theta)]=-\nabla R(\theta)$ for a differentiable objective function $R:\R^d \to \R$. As we discussed in the introduction, for noiseless optimization tasks, i.e. $X(U,\theta)\equiv -\nabla R(\theta)$, the \Adam\ vector field is given by
$$
f(\theta) = -\frac{\nabla R(\theta)}{|\nabla R(\theta)|+\eps}.
$$
In that case, it is easy to see that $R$ is a Lyapunov function for the set of critical points of $R$ under the semiflow \eqref{eq:ODE}, and one can use Corollary~\ref{cor:Lyapunov} to show that the limit set of the \Adam\ algorithm is a connected subset of the critical points. 

This is, in general, not the case when the innovation $(X,U)$ is random. In the stochastic case, $R$ might not be a Lyapunov function for the semiflow \eqref{eq:ODE} and the zero set of the \Adam\ vector field $f$ might be different from the set of critical points of $R$. 

Let us consider a stochastic setting motivated by empirical risk minimization in overparametrized supervised learning. In this case, the data points can be perfectly interpolated, so that the corresponding interpolating network achieves zero loss. Moreover, the variance of the gradient estimator with respect to the empirical distribution of the data points vanishes. As we will show in the next result, whenever the iterates of the \Adam\ algorithm are sufficiently close to such a minimizer, they are typically attracted to the set of minimizers, and hence the \Adam\ optimizer converges to the set of global optima.

More generally, we assume that there exists a set of minima for which the stochastic noise vanishes as the \Adam\ iterates approach this set. We prove that, whenever the noise decays sufficiently fast, $R$ serves as a Lyapunov function in a neighborhood of the minima, which in turn allows us to establish the convergence of the \Adam\ algorithm to this set.

\begin{theorem} \label{thm:degeneratenoiseNEW}
	Consider the setup of \cref{thm:item_ii} in \cref{thm:main}. Let $ p_1,p_2,q \in (0,\infty)$ with $1<p_2\le  2$,  $\sfrac{1}{2}p_2+1<p_1\le \frac 32 p_2$ and $q<\sqrt{\beta}$. Let $(X,U)$ be an innovation and 
	$R: \R^d \to \R$ be a $C^1$-function such that $\E[X(U, \theta)]=-\nabla R(\theta)$ for all $\theta \in \R^d$. 
	Let $\mathcal C$ be a closed, non-empty set of critical points and $K \subset \R^d $ be a compact set. Assume that there exists an open neighborhood $\mathcal D \supset \mathcal C \cap K$ and constants $C_{1}, C_{2}\in (0,\infty)$ such that 
	\bas{\label{eq:783256}\mathcal C \cap \overline{\mathcal D} = \arg\min R\big|_{\overline{\mathcal D}}=\{\theta \in \overline{\mathcal D}: \nabla R(\theta)=0\}} and
	for all $\theta \in \mathcal D$ and $i=1,\dots, d$
	\begin{align} \label{eq:998262}
		\E[|X^{(i)}(U, \theta)|^3] \le C_{1} |\nabla R(\theta)|^{p_1} \quad \text{ and } \quad \P(|X^{(i)}(U,\theta)|^2< C_{2} |\nabla R^{(i)}(\theta)|^{p_2})\le   q.
	\end{align}
	Then, there exists a neighborhood $V$ of $\mathcal C \cap K$ such that for almost all
		$$
		\omega \in \bigl\{\theta_n \in V  \text{ for infinitely many } n \in \N_0\bigr\} \cap  \bigcup_{n_0\in\N_0}\bigcap_{n \ge n_0}\{\theta_n \in K\}
		$$
		one has $\lim_{n \to \infty } d(\theta_n(\omega),\mathcal C)=0$.
\end{theorem}

We will prove Theorem~\ref{thm:degeneratenoiseNEW} in Section~\ref{sec:proof2}.

\begin{remark}
	The second inequality in \eqref{eq:998262} implies that for all $\theta\in \mathcal D$ and $i=1,\dots,d$ 
	$$
		\P(|X^{(i)}(U,\theta)|^2 \ge C_2 |\nabla R^{(i)}(\theta)|^{p_2})\ge 1-q>0,
	$$
	so that
	$$
		\E[|X^{(i)}(U,\theta)|^3] \ge (1-q) C_2^{\frac 32 } |\nabla R^{(i)}(\theta)|^{\frac {3}{2} p_2}.
	$$
	Thus, on a set $\mathcal D$ containing critical points of $R$, \eqref{eq:998262} can only be satisfied if $p_1\le \frac 32 p_2$. In the following, we will use Theorem~\ref{thm:degeneratenoiseNEW} for the choice $p_1=3$ and $p_2=2$.
\end{remark}

Lastly, let us discuss the empirical risk minimization problem for a typical regression task using artificial neural networks. 

\begin{definition}[Empirical risk minimization]
	Given
	\begin{itemize}
		\item $(y_1,z_1),\dots,(y_N,z_N)\in \mathcal Y\times \R$ with $\mathcal Y$ being an arbitrary set and $N\in\N$  (training data), and
		\item $\mathfrak N:\R^d\times \mathcal Y\to\R$ a function that is $C^1$ in the first component for every fixed second component (realization function)
	\end{itemize}
	we call \begin{align} \label{eq:ERM}
		R(\theta)&:= \frac{1}{2N} \sum_{i=1}^N  |\mathfrak N^\theta (y_i)-z_i|^2 \qquad (\theta \in \R^d , i=1, \dots, N)
	\end{align}
	the \emph{empirical loss function for the training data $(y_1,z_1),\dots,(y_N,z_N)$} (not mentioning the function $\mathfrak N$ and the space $\mathcal Y$ in the definition).
	
	The \emph{empirical gradient estimator for the training data $(y_1,z_1),\dots,(y_N,z_N)$} is the innovation $(X,U)$ with $U$ being uniformly distributed on $\{1,\dots,N\}$ and 
	\bas{
		X(u,\theta)=-\nabla \mathfrak N^\theta (y_u) (\mathfrak N^\theta (y_u)-z_u)
	}
	and the \emph{mini-batch gradient estimator for the training data $(y_1,z_1),\dots,(y_N,z_N)$ of size $M\in\N$} is the innovation $(X_M, \mathbf U)$ with $\mathbf U$ being a vector of $M$ independent uniformly distributed random variables  on $\{1,\dots,N\}$ and 
	\bas{ \label{eq:ERMmini-batch}
		X_M(\mathbf u,\theta)=-\frac 1M\sum_{i=1}^M \nabla \mathfrak N^\theta (y_{\mathbf u_i}) (\mathfrak N^\theta (y_{\mathbf u_i})-z_{\mathbf u_i}).
	}
\end{definition}

In empirical risk minimization the function $\mathfrak N$ is typically the realization function of an ANN. Provided that in the training data the entries $y_1,\dots,y_N$ are pairwise distinct and the ANN  is sufficiently expressive, one typically can find 
a $\theta^* \in \R^d$ with $R(\theta^*)=0$, see e.g.~\cite[Lemma~27.3]{foucart2022mathematical}.

In the latter regime, we will show that for sufficiently large mini-batch size $M$, the \Adam\ optimizer converges to an optimum, if it finds a neighborhood of the global minima infinitely often and the loss function $R$ satisfies a Polyak-\L ojasiewicz condition, introduced in \cite{polyak1963gradient, lojasiewicz1963propriete}. For the default value $\beta = 0.999$, a direct computation shows that $M > 0.001 N$ is sufficient in order for the next result to hold. 

\begin{theorem} \label{thm:ERM} 
	Let $\epsilon \in (0,\infty), \beta\in(\frac 14,1), \alpha \in [0,\sqrt \beta)$.
	Let $R$ be an empirical loss function for the training data $(y_1,z_1), \dots, (y_N,z_N)$ with pairwise different inputs $y_1, \dots, y_N \in \mathcal Y$. Consider the \Adam\ algorithm with damping parameters $(\alpha, \beta, \epsilon)$ and mini-batch gradient estimator $(X_M, \mathbf U)$ defined in \eqref{eq:ERMmini-batch} and assume that
	$$
	M >\frac{\log(2\sqrt \beta -1)}{\log(1-\frac 1N)}.
	$$
	Let $\mathcal C \subset R^{-1}(\{0\})$ be a closed and non-empty set and $K \subset \R^d$ be a compact set. Assume that there exists a neighborhood $\mathcal D \supset \mathcal C \cap K$ such that $\mathcal C \cap \mathcal D = \arg\min R\big|_{{\mathcal D}}$ and $R$ satisfies the Polyak-\L ojasiewicz inequality on $\mathcal D$, i.e. there exists a $\mu \ge 0$ such that for all $\theta \in \mathcal D$
	$$
	2\mu R(\theta) \le |\nabla R(\theta)|^2.
	$$
	Consider the setup of \cref{thm:item_ii} in \cref{thm:main}.
	Then there exists a neighborhood $V \supset \mathcal C\cap K$ such that for almost all
	$$
	\omega \in \bigl\{\theta_n \in V \text{ for infinitely many } n \in \N_0\bigr\} \cap  \bigcup_{n_0\in\N_0}\bigcap_{n \ge n_0}\{\theta_n \in K\}.
	$$
	one has $\lim_{n \to \infty } d(\theta_n(\omega),\mathcal C)=0$.
\end{theorem}

We will prove Theorem~\ref{thm:degeneratenoiseNEW} in Section~\ref{sec:proof3}.

\section{Proof of \cref{thm:main}} \label{sec:proof}
In this section, we prove the first main result of this article, \cref{thm:main}. We follow the ideas in \cite{dereich2024convergence} and define several intermediate approximations to the \Adam\ algorithm \eqref{eq:Adamdefi}. 
Throughout this section, we assume without further mentioning the following: 
The damping parameters $(\alpha,\beta,\epsilon)$ are fixed as in \cref{thm:main}, i.e.,  $0<\alpha < \sqrt \beta < 1$ and $\epsilon >0$. Moreover, $(\theta_n)_{n\in\N_0}$ is an \Adam\ algorithm started in $\theta_0\in\R^d$ with innovation $(X,U)$ and monotonically decreasing step-sizes $(\gamma_n)_{n\in\N}$ satisfying $\lim_{n\to\infty} \gamma_n=0$ and $\sum_{i=1}^\infty \gamma_i = \infty$, and $f$ denotes the corresponding \Adam\ vector field defined in \eqref{eq:fvector}.

Let $(\varrho_k)_{k \in -\N_0}$ be given by
\begin{align} \label{def:varrho}
\varrho_k = (1-\alpha) \epsilon^{-1} \Bigl( \alpha^{-k}+\frac{1}{\sqrt{1-\alpha^2/\beta}} \beta^{-k/2} \Bigr),
\end{align}
for all $k \in -\N_0$,
and denote by $\ell_\varrho^d$ the space of all $\R^d$-valued sequences  $(x_k)_{k \in -\N_0}$ satisfying
$$
\|(x_k)_{k \in -\N_0}\|_{\ell_\varrho^d} := \sum_{k=-\infty}^0 \varrho_k |x_k| < \infty.
$$
We define  $g: \ell^d_\varrho \to \R^d$ via
$$
g^{(j)}((x_k)_{k \in -\N_0}) =\frac{(1-\alpha)\sum_{k=-\infty}^0\alpha^{-k}x_k^{(j)}}{\sqrt{(1-\beta)\sum_{k=-\infty}^0 \beta^{-k}(x_k^{(j)})^2}+\epsilon},
$$
for all $j=1, \dots, d$, and recall the following fact from \cite{dereich2024convergence} that we will frequently use throughout this article.

\begin{lemma} \label{lem:lipschitz}
	The mapping $g: \ell_\varrho^d \to \R^d$ is $1$-Lipschitz and uniformly bounded by $\sqrt d \frac{1-\alpha}{\sqrt{1-\beta}} \frac{1}{\sqrt{1-\alpha^2/\beta}}$. Moreover, if the innovation $(X, U)$ is $p$-regular for a $p \in [2, \infty)$, see Definition~\ref{assu:regular}, then the corresponding \Adam\ vector field $f: \R^d \to \R^d$ defined in Definition~\ref{def:Adam} is uniformly bounded by $\sqrt d \frac{1-\alpha}{\sqrt{1-\beta}}\frac{1}{\sqrt{1-\alpha^2/\beta}}$ and locally Lipschitz-continuous.
\end{lemma}

\begin{proof}
	See Lemma~3.1 and Equation (148) in \cite{dereich2024convergence}.
\end{proof}

For simplicity we write for $k \in \mathbb Z$
$$
X_k = \begin{cases}
	X(U_k,\theta_{k-1}), & \text{ if } k \ge 1, \\
	0, & \text{ if } k \le 0,
\end{cases}
$$
and $\textbf X(k)= (X_{k+i})_{i \in -\N_0}$.
We fix an increasing sequence of natural numbers $(n_\ell)_{\ell \in \N_0}$ and  define the processes
\begin{align} 
\label{eq:approximations}
  \Theta_n := \sum_{k=1}^n \gamma_k g(\textbf X(k)) 
  \quad, 
  \quad 
  \tilde \Theta_n := \sum_{k=1}^n \gamma_k g(\tilde {\textbf X}(k)) 
  \quad, 
  \quad \dbtilde \Theta_n := \sum_{k=1}^n \gamma_k g( \dbtilde {\bf X}( k ) ) ,
\end{align}
where for $n \in \{n_{\ell-1}+1, \dots, n_\ell\}$ we define $\tilde {\textbf X}(n)= (\tilde {\textbf X}(n)_{k})_{k \in -\N_0}$ and $\dbtilde {\textbf X}(n)= (\dbtilde {\textbf X}(n)_{k})_{k \in -\N_0}$ via 
$$
\tilde {\textbf{X}}(n)_k = \begin{cases} 
	X(U_{n+k}, \theta_{n_{\ell-1}}), & \text{ if } n+k > n_{\ell-1}, \\
	X_{n+k}, & \text{ else},
\end{cases}
$$
and
\begin{align} \label{eq:dbtildeX}
	\dbtilde {\textbf{X}}(n)_k = \begin{cases} 
		X(U_{n+k}, \theta_{n_{\ell-1}}), & \text{ if } n+k > n_{\ell-1}, \\
		X(\tilde U_{n+k}, \theta_{n_{\ell-1}}), & \text{ else}.
	\end{cases}
\end{align}
Here, $(\tilde U_{n})_{n \in \mathbb Z}$ denotes a family of independent copies of $U$ that is independent of $(U_n)_{n \in \N}$. 

The processes $(\Theta_n)_{n \in \N_0}, (\tilde \Theta_n)_{n \in \N_0}$ and $(\dbtilde \Theta_n)_{n \in \N_0}$ can be viewed as approximations to the \Adam\ algorithm \eqref{eq:Adamdefi}.  In fact, the increments $\Delta \theta_n:= \theta_n-\theta_{n-1}$ and $\Delta \Theta_n :=\Theta_n-\Theta_{n-1}$ only differ in the non-presence of the bias correction term $1-\beta^n$ in the definition of $\sigma_n$.  The process $(\tilde \Theta_n)_{n \in \N_0}$ evaluates the innovation $(X,U)$ within the current time frame $(n_{\ell-1}, n_\ell]$ always with the $\theta$-parameter being frozen at the beginning of the time frame. It can be conceived as Euler approximation to $(\Theta_n)_{n \in \N_0}$ with the updates occurring at the time instances  in $(n_\ell)_{\ell \in \N_0}$.
The process $(\dbtilde \Theta_n)_{n \in \N_0}$ applies an additional randomization and replaces the true historical information outside the current time frame by an independently generated history. We will carry out an error analysis and argue that the behavior of $(\theta_n)_{n \in \N_0}$ is closely related to the behavior of $(\dbtilde \Theta_n)_{n \in \N_0}$. Afterwards, we will use that a law of large numbers applies for the action of the process   $(\tilde{\tilde \Theta}_n)_{n \in \N_0}$ on the time frames. For this note that $\theta_{n_\ell}$ is approximately equal to $\theta_{n_{\ell-1}} + \dbtilde\Theta_{n_\ell}-\dbtilde\Theta_{n_{\ell-1}}$ and we have that
\bas{ \label{eq:1121}
	\theta_{n_{\ell-1}} + \E\bigl[\dbtilde\Theta_{n_\ell}- \dbtilde\Theta_{n_{\ell-1}} | \theta_{n_{\ell-1}} \bigr]= \theta_{n_{\ell-1}} + (t_{n_{\ell}} -t_{n_{\ell-1}} ) \,f( \theta_{n_{\ell-1}})
}
which is close to $\Psi_{t_{n_{\ell}} -t_{n_{\ell-1}} }( \theta_{n_{\ell-1}} )$ as long as the time increment is small.

The proof of \cref{thm:main} is divided into the following steps. 
First, we quantify the individual approximation errors using Proposition~5.2 in \cite{dereich2024convergence}. In Section~\ref{sec:pregular}, we use the local $p$-regularity of the innovation $(X,U)$, to bound the influence of the past gradient evaluations. In Section~\ref{sec:summability}, we use the assumptions on the step-sizes to show a summability result for the approximation errors. Finally, in Section~\ref{sec:proofmain} we finish the proof using a Borel-Cantelli argument.

Let us present a slight variation of \cite[Proposition 5.2]{dereich2024convergence} which quantifies the individual approximation errors of the intermediate approximations.

\begin{proposition}[Proposition~5.2 in \cite{dereich2024convergence}]\label{prop:234}
	Let $C, L\in[0,\infty)$, $p\in[2,\infty)$ and $V\subset\R^{d}$ a measurable set such that for every $\theta,\theta'\in V$, one has
	$$
	\E[|X(U,\theta)|^{p}]^{1/p}\le C,  \ \text{ \ and}
	\qquad \E[|X(U,\theta)-X(U,\theta')|^{p}]^{1/p} \le  L  \,|\theta-\theta'|.
	$$
	Moreover, let
		$$
			\mathfrak N=\inf \{n\ge n_{0}:  \theta(n)\not \in V\}
		$$
		One has the following for the approximations $(\Theta_n)_{n\ge n_0}$, $(\tilde\Theta_n)_{n\ge n_0}$ and $(\dbtilde\Theta_n)_{n\ge n_0}$:
		\begin{enumerate}
			\item[(I)] for every $\mathfrak n\in\{n_{0},n_{0}+1,\dots\}$,
			\begin{align}\begin{split}
					\E\Bigl[ \Bigl(\sum_{k=\mathfrak n+1}^{\infty} \1_{\{\mathfrak N\ge k\}} |\Delta \theta_{k}-\Delta \Theta_{k} |\Bigr)^{p}\Bigr]^{1/p}&\le \kappa_1 \gamma_{\mathfrak n+1} \beta^{\mathfrak n+1}
			\end{split}
		\end{align}
			\item[(II)] for every $\ell\in\N$ and $\mathfrak n,n\in\N_{0}$ with $n_{\ell-1}\le \mathfrak n\le n\le n_{\ell}$,
			\begin{align}\begin{split}
					\E\Bigl[ \Bigl(\sum_{k=\mathfrak n+1}^n \1_{\{\mathfrak N\ge k\}} |\Delta \Theta_{k} -\Delta \tilde\Theta_{k}|\Bigr)^{p}\Bigr]^{1/p}&\le \kappa_{2} L  (t_{n_{\ell}}-t_{n_{\ell-1}})(t_{n}-t_{\mathfrak n})
			\end{split}\end{align}
			\item[(III)] for every $\ell\in\N$,
			\begin{align}
				\label{eq78246-2}
				\E\Bigl[&\1_{\{\mathfrak N>n_{\ell-1}\}}\Bigl(\sum _{k=n_{\ell-1}+1}^{n_{\ell}} \bigl|\Delta {\tilde \Theta}_{k}-\Delta \dbtilde\Theta_{k}\bigr|\Bigr)^{p} \Big| \mathcal F_{n_0}\Bigr]^{1/p}\le \gamma_{n_{\ell-1}+1} \bigl(\kappa_{3} C+\kappa_{4}\beta^{(n_{\ell-1}-n_{0})/2}\|\mathbf X(n_0)\|_{\ell_{\varrho}^{d}}),
			\end{align}
			\item[(IV)]for every $\ell\in\N$,
			$$
				\E[\1_{\{\mathfrak N>n_{\ell-1}\}}|\dbtilde\Theta_{n_{\ell}}-\dbtilde \Theta_{n_{\ell-1}}-(t_{n_{\ell}}-t_{n_{\ell-1}})f(\theta_{n_{\ell-1}})|^{p}|\cF_{n_{\ell-1}}]^{1/p} \le \kappa_{5} C\Bigl(\sum_{k=n_{\ell-1}+1}^{n_{\ell}} \gamma_k^{2}\Bigr)^{1/2},
			$$
			\end{enumerate}
			where 	
		\bas{ \label{eq237643}
			&\kappa_{1}=\frac2{1-\beta}\sqrt d \frac{1-\alpha}{\sqrt{1-\beta}} \frac{1}{\sqrt{1-\alpha^2/\beta}}, \  \kappa_{2}=\frac{1-\alpha}{\sqrt{1-\beta}\sqrt{1-\alpha^{2}/\beta}} d \|\varrho\|_{\ell_{1}},\\ 
			&\kappa_{3}=2 \frac{\sqrt\beta}{1-\sqrt\beta}\|\varrho\|_{\ell_{1}}, \  
			\kappa_{4}=\frac{\sqrt \beta}{1-\sqrt\beta},  \ \kappa_{5}=4  C_{p}\bigl((1-\sqrt\beta)^{-1} \sqrt{\varrho_{0} \|\varrho\|_{\ell_{1}}}+ \|\varrho\|_{\ell_{1}} \bigr),}
		$(\varrho_{k})_{k \in -\N_0}$ is as in~\eqref{def:varrho} and $C_{p}$ is the constant in the Burholder-Davis-Gundy inequality when applied for the $p$th moment.
	\end{proposition}
	
	\begin{proof}
		We make only two minor modifications to Proposition~5.2 in \cite{dereich2024convergence}. To get (I) we use boundedness of $g$ and to get (IV) we use the fact that
		\begin{equation}
			\E[|X(U,\theta)-E[X(U,\theta)]|^p]^{1/p} \le \E[|X(U,\theta)|^p]^{1/p} + |\E[X(U,\theta)]| \le 2\E[|X(U,\theta)|^p]^{1/p}.
		\end{equation}
	\end{proof}

\subsection{Bounding the influence from past evaluations using the growth bound} \label{sec:pregular}

In order to apply Proposition~\ref{prop:234}, especially \eqref{eq78246-2}, in the proof of the main result (\cref{thm:main}), we establish a bound for the $p$th moment of $\|\mathbf X(n_0)\|_{\ell_{\varrho}^{d}}$. The next lemma makes use of the growth bound for the innovation $(X,U)$ given in \eqref{eq:34296255544} to control the effect of past simulations collected prior to the comparison horizon. Exploiting the boundedness of the velocity of the \Adam\ algorithm, we arrive at the following statement.

\begin{lemma} \label{lem:growthbound}
	Let $p \in[ 1,\infty)$ and
	suppose	that  there exists a constant $C_1\in(0,\infty)$ so that for every $\theta\in\R^d$
	\[
	\E[|X(\theta,U)|^p]^{1/p}  \le  C_1(e^{C_1|\theta|}+1).
	\]
	Then,  for every $R\in(0,\infty)$, one has that
	\[
	\sup_{N\in\N_0}
	\E[\1_{\{|\theta_N|\le R\}}\|\mathbf X(N)\|_{\ell_\varrho}^p]<\infty.
	\]
\end{lemma}

\begin{proof}
	Due to the uniform boundedness of $g: \ell_\varrho^d \to \R^d$ there exists a constant $\kappa_1\ge 0$ such that for all $n \in \N$
	\begin{align} \label{eq:238942936822}
		|\theta_n-\theta_{n-1}|\le \kappa _1\gamma_n.
	\end{align}
	Consequently, on the event $\{|\theta_N|\le R\}$ we get that for every $k\in\{0,\dots,N\}$
	\[
	|\theta_k|\le R+ \kappa_1 (t_N-t_k)\le R+ \kappa_1  \frac{t_N}N (N-k),
	\]
	where we have used that $t_k \ge \frac{k}{N}t_N$ due to the monotonicity of $(\gamma_n)_{n \in \N}$.
	Thus,
	\begin{align*}
		\E\Big[\1_{\{|\theta_N|\le R\}}|X(\theta_{k-1},U_k)|^p\Big]^{1/p}&\le \E\Big[\1_{\{|\theta_{k-1}|\le R+ \kappa_1  \frac{t_N}N (N-k+1) \}} |X(\theta_{k-1},U_k)|^p\Big]^{1/p}\\
		&\le C_1 \Big(\exp\Big\{C_1\Big(R+ \kappa _1 \frac{t_N}N (N-k+1)\Big)\Big\}+1\Big).
	\end{align*}
	Next, note that since $\alpha < \sqrt \beta$ there exists a $\kappa_2\in(0,\infty)$ so that for every $\mathbf x \in \ell_\varrho^d$
	\[\|\mathbf x\|_{\ell_\varrho^d} \le \kappa_2 \sum_{k\in-\N_0} \beta^{-k/2} |x_k|.
	\]
	Using the triangle inequality for the $L^p$-norm,
	\bas{\label{eq87346}
		&\E[\1_{\{|\theta_N|\le R\}}\|\mathbf X_N\|_{\ell_\varrho^d}^p]^{1/p}
		\le  \E\Big[\1_{\{|\theta_N|\le R\}}\kappa_2^p \Big( \sum_{k=1}^N \beta^{(N-k)/2} |X(\theta_{k-1},U_k)|\Big)^p\Big]^{1/p} \\		
		&\le  C_1  \kappa_2\sum_{k=1}^N \beta^{(N-k)/2} \Big(\exp\Big\{C_1\Big(R+ \kappa _1 \frac{t_N}N (N-k+1)\Big)\Big\}+1\Big)\\
		&= C_1 \kappa_2 \sum_{k=1}^N  \Bigl(\exp\Big\{-\Big(\frac 12\log \beta^{-1}-C_1 \kappa_1\frac{t_N}N\Big)(N-k)+ C_1\Big(R+ \kappa _1 \frac{t_N}N\Big)\Big\}+\beta^{(N-k)/2}\Big)
	}
	Since $\gamma_n \to 0$ one has that $t_N/N \to 0$.
	Thus, there exist $N_0\in\N_0$ and $\delta,\kappa_3\in(0,\infty)$ such that for all $N\ge N_0$
	\[
	-\Big(\frac 12\log \beta^{-1}-C_1 \kappa_1\frac{t_N}N\Big)\le -\delta \text{ \ and \  } C_1\Big(R+ \kappa _1 \frac{t_N}N\Big)\le \kappa_3,
	\]
	which implies that for those $N$
	\[
	\E[\1_{\{|\theta_N|\le R\}} \|\mathbf X_N\|_{\ell_\varrho}^p]^{1/p}
	\le C_1 \kappa_2 \sum_{k\in\N_0} (e^{\kappa_3-\delta k}+
	\beta^{k/2})<\infty.
	\]
	For the remaining finitely many $N$'s the expectations are finite as consequence of~(\ref{eq87346}).
\end{proof}

\subsection{Summability conditions} \label{sec:summability}
In this section, we prove summability of certain sequences that appear as error bounds of the intermediate approximations using the assumptions in \cref{thm:main}. First, we prove a consequence of the step-size condition in \cref{thm:main}. 
\begin{lemma}\label{lem:newgamma}
	Let $(\gamma_n)_{n \in \N_0}$ be a  non-summable decreasing sequence of reals with limit zero and suppose that
	\bas{ \label{eq:209346237522}
		\exists\varepsilon_1,\varepsilon_2\in(0,\infty), n_0\in\N \,\forall n\ge n_0: \gamma_{n+\lceil\varepsilon_1\gamma_n^{-1}\rceil}\ge \varepsilon_2 \gamma_n.
	}
	Let $\Gamma: (0,\infty) \to [0,\infty)$ be defined by setting $\Gamma_{t} = \gamma_n$ for all $n \in \N$ and $t \in (t_{n-1},t_{n}]$, where $t_n=\sum_{i=1}^n \gamma_i$. 
	Then one has for 
	every $c\in(0,\infty)$ and $q \in (0,1)$
	$$
		\limsup_{t\to\infty} \frac {\Gamma_t}{\Gamma_{t+c}}<\infty
		\text{ \ \ and \ \ }\limsup_{n\to\infty} \frac{\gamma_n}{ \gamma_{n+\lceil c \gamma_n^{-q}\rceil}}<\infty .
	$$
\end{lemma}

\begin{proof}The second statement is an immediate consequence of the assumption and we only prove the first statement.
	Fix $\varepsilon_1,\varepsilon_2\in(0,\infty)$, $n_0\in\N$ such that  for all $n\ge n_0$, 
	\bas{\label{eq:34476356} \gamma_{n+\lceil\varepsilon_1\gamma_n^{-1}\rceil}\ge \varepsilon_2 \gamma_n.}
	Using monotonicity of $(\gamma_n)_{n \ge n_0}$ we get that for  every $n\ge n_0$
	\bas{\label{eq:72452}
		t_{n+\lceil\varepsilon_1\gamma_n^{-1}\rceil}\ge t_n + \lceil\varepsilon_1\gamma_n^{-1}\rceil 
		\,\gamma_{n+\lceil\varepsilon\gamma_n^{-1}\rceil}\ge t_n+ \lceil\varepsilon_1\gamma_n^{-1}\rceil 
		\, \varepsilon_2 \gamma_n\ge t_n+\varepsilon_1\varepsilon_2.
	}
	Now pick $\varepsilon\in(0,\varepsilon_1 \varepsilon_2)$ and $n_1\in\N$ with $n_1\ge n_0$ and $\gamma_{n_1}\le \varepsilon_1 \varepsilon_2-\varepsilon$. Fix  $t\in[t_{n_1},\infty)$ and $n\in \N$ with $t\in(t_{n-1},t_n]$. Then $n\ge n_1$ and we have $t_n+\varepsilon_1\varepsilon_2 \ge t+\varepsilon$ so that
	$$
		\Gamma_{t+\varepsilon}\ge \Gamma_{t_n+\varepsilon_1\varepsilon_2}\ge \gamma_{n+\lceil\varepsilon_1\gamma_n^{-1}\rceil}\ge \eps_2\gamma_n =\eps_2 \Gamma_t,
	$$
	where we have used (\ref{eq:72452}) in the second inequality and (\ref{eq:34476356}) in the third.
	Now choose, for fixed  $c\in(0,\infty)$, an $m\in\N$ with $m\eps\ge c$ and note that for all $t\ge t_{n_1}$
	$$
		\frac {\Gamma_t}{\Gamma_{t+c}}\le \frac {\Gamma_t}{\Gamma_{t+m\eps}} = \frac {\Gamma_t}{\Gamma_{t+\eps}}\ldots  \frac {\Gamma_{t+(m-1)\eps}}{\Gamma_{t+m\eps}}\le \eps_2^{-m}.
	$$
\end{proof}

The definition of our approximations in (\ref{eq:approximations}) make use of an  increasing sequence $(n_\ell)_{\ell \in \N_0}$ of natural numbers. To keep the error of the approximations balanced we work with particular sequences to be defined next.

\begin{definition} \label{def:partition}
	Let $n_0 \in \N_0$ and $\rho \in [\gamma_{n_0+1}^{2/3},\infty)$. We call the $\N_0$-valued sequence $(n_\ell)_{\ell \in \N_0}$ a \emph{$\rho$-partition} if, for every $\ell \in \N$, $n_\ell$ is the largest integer that satisfies 
	$$
	t_{n_\ell}-t_{n_{\ell-1}}\le \rho \gamma_{n_{\ell-1}+1}^{1/3}.
	$$
\end{definition}

The proof of statement \ref{thm:item_ii}  of \cref{thm:main} is based on a Borel-Cantelli argument. The necessary summability condition will be verified using the following result. 

\begin{proposition} \label{prop:summable}
	Let $p \in[ 1,\infty)$ and $(\gamma_n)_{n\in\N}$ be a non-summable, decreasing sequence of reals with limit zero such that $\sum_{n\in \N} \gamma_n^{1+\frac p3}<\infty$ and \eqref{eq:209346237522} holds. 	Then there exists a $\rho \in(0,\infty)$ and a $\rho$-partition $(n_\ell)_{ \ell \in \N_0}$ such that the following holds:
	
	Denote by  $\II_s=(t_{n_{\ell(s)}},t_{n_{r(s)}}]$, $s\in\N$, disjoint intervals of length at most $2T$.
	Then the quantities
	\bas{ \label{eq:mathcalE}
		\cE_2(s) =  &\sum_{j=\ell(s)+1}^{r(s)} (t_{n_{j}}-t_{n_{j-1}})^2 \ , \ 
		\mathcal E_3(s)= \sum_{j=\ell(s)+1}^{r(s)} \gamma_{n_{j-1}+1} \\
		& \ \text{ and } \  \mathcal E_4(s) = \sum_{j=\ell(s)+1}^{r(s)} \Bigl(\sum_{k=n_{j-1}+1}^{n_{j}} \gamma_k^{2}\Bigr)^{1/2}
	}
	are in $\ell_p$ meaning that
	$$
	\sum_{s\in\N} \bigl(\cE_2(s)^p+\cE_3(s)^p+\cE_4(s)^p\bigr) < \infty.
	$$
\end{proposition}

\begin{proof} 
	By Lemma~\ref{lem:newgamma}, there exists a $\kappa_1 >0$ such that for all $n \in \N_0$
	$$
	\gamma_{n+ \lceil \gamma_{n+1}^{-2/3}\rceil}\ge  \kappa_1\gamma_{n+1}.
	$$
	Let $\rho < \kappa_1$ and $(n_{\ell})_{\ell \in \N_0}$ be a $\rho$-partition. Thus, for all $\ell \in \N$, we have that
	$$
	t_{n_{\ell-1} + \lceil \gamma_{n_{\ell-1}+1}^{-2/3}\rceil} -t_{n_{\ell-1}} \ge 
	\kappa_1 \gamma_{n_{\ell-1}+1} \lceil \gamma_{n_{\ell-1}+1}^{-2/3}\rceil.
	$$
	so 
	that 
	\bas{ \label{eq:983730003}
		n_{\ell}\le n_{\ell-1}+  \gamma_{n_{\ell-1}+1}^{-2/3}\text{ \ and \ } \gamma_{n_\ell}\ge \gamma_{n_{\ell-1}+ \lceil  \gamma_{ n_{\ell-1}+1}^{-2/3}\rceil}\ge \kappa_1 \gamma_{n_{\ell-1}+1}.
	}

	By definition of the $\rho$-partition, the summands in the definition of $\cE_2$ satisfy $(t_{n_j}-t_{n_{j-1}})^2\sim \rho^2 \gamma_{n_{j-1}+1}^{2/3}$ which is larger than $\gamma_{n_{j-1}}$ for all sufficiently large $j\in\N$. So for sufficiently large $s\in\N$, we get that  $\cE_3(s)\le \cE_2(s)$ and it suffices to prove finiteness of $\sum_{s\in\N}\mathcal E_2(s)^p$ and $\sum_{s\in\N}\mathcal E_4(s)^p$.

	Using Jensen's inequality and \eqref{eq:983730003} we get
	\bas{ \label{eq:238969000}
		\cE_2(s) & \le \sum_{j=\ell(s)+1}^{r(s)} \Big(\sum_{i=n_{j-1}+1}^{n_{j}}\gamma_i\Big)^2 =  \sum_{j=\ell(s)+1}^{r(s)} (n_{j}-n_{j-1})^2 \Big(\frac{1}{n_{j}-n_{j-1}} \sum_{i=n_{j-1}+1}^{n_{j}}\gamma_i\Big)^2 \\
		&\le \sum_{j=\ell(s)+1}^{r(s)} (n_{j}-n_{j-1})  \sum_{i=n_{j-1}+1}^{n_{j}}\gamma_i^2 \\
		& \le \sum_{k=n_{\ell(s)}+1}^{n_{{r(s)}}} \gamma_k^{4/3}.
	}
	We denote by $N_s=n_{{r(s)}}-n_{{\ell(s)}+1}$ the number of summands in the latter sum and and, again, apply Jensen's inequality to arrive at
	$$
	\cE_2(s)^p \le N_s^p\Bigl(\frac 1{N_s} \sum_{k=n_{\ell(s)}+1}^{n_{{r(s)}}}  \gamma_k^{4/3}\Bigr)^p\le  N_s^{p-1}  \sum_{k=n_{\ell(s)}+1}^{n_{{r(s)}}}\gamma_k^{\frac 43p} 
	$$
	By Lemma~\ref{lem:newgamma}, there exists a constant $\kappa_2>0$ such that for all $s\in\N$ 
	$$
	\gamma_{n_{\ell(s)+1}+1} \le \kappa_2 \gamma_{n_{r(s)}}.
	$$
	Together with  $N_s \gamma_{n_{r(s)}}\le 2T$ we obtain that $N_s\le 2\kappa_2 T/\gamma_{n_{\ell(s)}+1}$ so that
	$$
	\cE_2(s)^p \le   (2\kappa_2 T)^{p-1} \gamma_{n_{\ell(s)}+1}^{-(p-1)}  \sum_{k=n_{\ell(s)}+1}^{n_{{r(s)}}}\gamma_k^{\frac 43p} \le (2\kappa_2 T)^{p-1}   \sum_{k=n_{\ell(s)}+1}^{n_{{r(s)}}}\gamma_k^{1+\frac p3}.
	$$
	Since every summand $k$ occurs in at most one term of $(\cE_2(s))_{s \in \N}$ we conclude that
	$$
	\sum_{s\in\N}\cE_2(s)^p \le (2 \kappa_2 T)^{p-1}  \sum_{k\in\N} \gamma_k^{1+\frac p3} < \infty.
	$$
	
	In order to prove $\sum_{s\in\N}\mathcal E_4(s)^p<\infty$,
	first note that 
	$$
		\Bigl(\sum_{k=n_{j-1}+1}^{n_{j}} \gamma_k^{2}\Bigr)^{1/2} \le \gamma_{n_{j-1}+1} \sqrt {n_{j}-n_{j-1}}.
	$$
	Thus, we can use the first inequality in \eqref{eq:983730003} to get
	\begin{align*}
	\mathcal E_4(s) &\le  \sum_{j=\ell(s)+1}^{r(s)} \gamma_{n_{j-1}+1}^{2/3} = \sum_{j=\ell(s)+1}^{r(s)} \gamma_{n_{j-1}+1}^{2/3} \Biggl( (t_{n_j}-t_{n_{j-1}})^{-1} \sum_{k=n_{j-1}+1}^{n_j} \gamma_k\Biggr).
	\end{align*}
	Furthermore,
	$$
		t_{n_j}-t_{n_{j-1}}+\gamma_{n_j+1} \ge \rho \gamma_{n_{j-1}+1}^{1/3}
	$$
	so that for sufficiently large $j \in \N$
	$$
		t_{n_j}-t_{n_{j-1}}\ge \frac{1}{2} \rho \gamma_{n_{j-1}+1}^{1/3}.
	$$
	Consequently, we have for sufficiently large $s$ (say for $s\ge s_0$) that
	\begin{align} \label{eq:23975843}
		\mathcal E_4(s) \le \sum_{j=\ell(s)+1}^{r(s)} 2 \rho^{-1} \gamma_{n_{j-1}+1}^{1/3} \Biggl(\sum_{k=n_{j-1}+1}^{n_j} \gamma_k\Biggr) \le  \frac 2{\rho \kappa_1^{1/3}}  \sum_{k=n_{\ell(s)}+1}^{n_{{r(s)}}} \gamma_k^{4/3}.
	\end{align}
	Up to constants this estimate is the same as the one that we obtained in (\ref{eq:238969000}) for $\cE_2(s)$. Based on this inequality we proved finiteness of the $\cE_2$-series and for the same reason the $\cE_4$-series converges to a finite value.
\end{proof}

\subsection{Proof of \cref{thm:main}} \label{sec:proofmain}
	We are now in the position to prove \cref{thm:main}. We first prove statement  \ref{thm:item_ii}.
	
	\underline{Step 1:} Using the intermediate approximations.
	
	Fix $T, R\ge 0$. Let $\rho >0$ and $(n_\ell)_{\ell \in \N_0}$ be a $\rho$-partition such that the conclusion of Proposition~\ref{prop:summable} is true. By Definition~\ref{def:partition}, we have
	\begin{align} \label{eq:234675623}
		\sum_{i=n_{\ell-1}+1}^{n_\ell} \gamma_i = t_{n_\ell}-t_{n_{\ell-1}} \le \rho \gamma_{n_{\ell-1}+1}^{1/3} \overset{\ell \to \infty}{\longrightarrow} 0.
	\end{align}
	For $K \in \N$ we define $N_K := \inf\{\ell \in \N: \sum_{i=1}^{n_\ell} \gamma_i \ge KT  \}$. For all but finitely many $K$'s one has $T/2 \le N_{K+1}-N_{K}\le 2T$ and, without loss of generality we assume that this is true for all $K\in \N$.
	First, we prove that it suffices to show that almost surely
	\begin{align} \label{eq:2834656252}
		\limsup_{K \to \infty} \, \1_{\{|\theta_{n_{N_{K-1}}}|\le R\}} \sup_{\ell=N_{K-1}, \dots, N_{K}-1} \sup_{i=\ell+1, \dots, N_{K}}  |\theta_{n_i} - \Psi_{t_{n_i}-t_{n_\ell}}(\theta_{n_\ell})| =0.
	\end{align}
	
	Recall that, by Lemma~\ref{lem:dereich2}, there exists a constant $C_1\ge 0$ with $|\theta_m-\theta_n|\le C_1 (t_m-t_n)$, almost surely, for all $0\le n\le m$ and $|\Psi_t(\theta)-\theta|\le C_1t$ for all $t\ge0$ and $\theta \in \R^d$, as well as a constant $C_2\ge 0$ such that for all $\theta,\theta' \in \overline {B_{R+4C_1T}(0)}$ one has $|f(\theta)-f(\theta')|\le C_2 |\theta-\theta'|$. Thus, for all $i \in \N$ and $k \in \{n_{i}+1, \dots, n_{i+1}\}$ one has
	$$
		|\theta_k-\theta_{n_i}| \le C_1( t_{n_{i+1}}-t_{n_i} ) \to 0
	$$
	and, moreover, for $K \in \N$, $i \in N_{K-1}, \dots, N_K-1$, $k \in \{n_i+1, \dots, n_{i+1}\}$, on the event $\{|\theta_{n_{N_{K-1}}}|\le R\}$, one has
	$
	|\theta_{n_i}|\le R+2C_1T
	$
	and
	$
	|\theta_k|\le R+2C_1T,
	$
	so that for all $s \in [t_{n_{i+1}}-t_k,t_{n_{i+1}}-t_{n_i}]$
	\begin{align*}
		|\Psi_{s}(\theta_{n_i}) - \Psi_{s-t_k+t_{n_i}}(\theta_k)| &\le \int_0^{s-t_k+t_{n_i}} |f(\Psi_{u}(\Psi_{t_k-t_{n_{i}}}(\theta_{n_i})))-f( \Psi_{u}(\theta_k))| \, du \\
		& \le \int_0^{s-t_k+t_{n_i}} C_2 |\Psi_{u}(\Psi_{t_k-t_{n_{i}}}(\theta_{n_i}))- \Psi_{u}(\theta_k)| \, du.
	\end{align*}
	Thus, using Gronwall's inequality,
	$$
		|\Psi_{s}(\theta_{n_i}) - \Psi_{s-t_k+t_{n_i}}(\theta_k)| \le C_1 (t_{n_{i+1}}-t_{n_i}) \exp(C_2 (t_{n_{i+1}}-t_{n_i})) \to 0.
	$$
	We deduce that 
	\begin{align*} 
		&\limsup_{K \to \infty} \,\1_{\{|\theta_{n_{N_{K-1}}}|\le R\}}  \; \sup_{\ell=N_{K-1}, \dots, N_{K}-1} \sup_{i=\ell+1, \dots, N_{K}}  |\theta_{n_i} - \Psi_{t_{n_i}-t_{n_\ell}}(\theta_{n_\ell})| \\
		& = \limsup_{K \to \infty} \, \1_{\{|\theta_{n_{N_{K-1}}}| \le R\}} \sup_{m=n_{N_{K-1}}, \dots, n_{N_{K}}-1} \; \sup_{k=m+1, \dots, n_{N_{K}}}  |\theta_{k} - \Psi_{t_{k}-t_{m}}(\theta_{m})|.
	\end{align*}
	Now, if $\{n_{N_{K-1}}, \dots, n_{N_{K}}-1\} \neq \emptyset$ then for every $N \in \{n_{N_{K-1}}, \dots, n_{N_{K}}-1\}$ it holds that 
	$$
	\{n \ge N : t_n-t_N \le T\} \subset [n_{N_{K-1}}, n_{N_{K+2}}]
	$$
	so that
	\begin{align*}
		&\limsup_{N \rightarrow \infty} \, \1_{\{|\theta_N| \le R\}} \sup_{\substack{n\in\N_0: \\ 0 \le t_n-t_N \leq T}}\left| \theta_{n}-\Psi_{t_n-t_N}( \theta_{N})\right| \\
		& \le 3 \limsup_{K \to \infty} \, \1_{\{|\theta_{n_{N_{K-1}}}| \le R+2C_1T\}} \sup_{m=n_{N_{K-1}}, \dots, n_{N_{K}}-1} \; \sup_{k=m+1, \dots, n_{N_{K}}}  |\theta_{k} - \Psi_{t_{k}-t_{m}}(\theta_{m})|.
	\end{align*}
	This proves the claim.

	In order to prove \eqref{eq:2834656252}, let us fix a $K \in \N$. Again, $|\theta_{n_{N_{K-1}}}| \le R$ implies $|\theta_{n_\ell}|\le R+2C_1T$ for all $\ell \in N_{K-1}, \dots, N_K$. 
	For all $i,\ell \in \{N_{K-1}, \dots, N_K\}$ with $i > \ell$ one has
	\begin{align*}
		|\theta_{n_i} - \Psi_{t_{n_i}-t_{n_\ell}}(\theta_{n_\ell})|&\le C \,  \sum_{j = \ell+1}^{i} \Big| E_1(j)+E_2(j)+E_3(j)+E_4(j)+E_5(j) \Big|,
	\end{align*}
	where
	\begin{align}
		\begin{split} \label{eq:2930472042}
			& E_1(j) = \sum_{k=n_{j-1}+1}^{n_j} \Delta \theta_k-\Delta \Theta_k \quad , \quad 
			E_2(j) = \sum_{k=n_{j-1}+1}^{n_j} \Delta \Theta_k-\Delta \tilde \Theta_k \\
			& E_3(j) = \sum_{k=n_{j-1}+1}^{n_j} \Delta \tilde \Theta_k-\Delta \dbtilde \Theta_k \quad , \quad 
			E_4(j) = \dbtilde\Theta_{n_{j}}-\dbtilde \Theta_{n_{\ell-1}}-(t_{n_{j}}-t_{n_{j-1}})f(\theta_{n_{j-1}}) \\
			& \text{and } \qquad E_5(j) = (t_{n_j}-t_{n_{j-1}})f(\theta_{n_{j-1}}) -  (\Psi_{t_{n_j}-t_{n_{j-1}}}(\theta_{n_{j-1}}) - \theta_{n_{j-1}}).
		\end{split}
	\end{align}
	Using Proposition~\ref{prop:234} and Lemma~\ref{lem:growthbound},
	there exists $\kappa_1, \dots, \kappa_5$ such that
	\begin{align*}
		\E \Big[\1_{\{|\theta_{n_{N_{K-1}}}| \le R\}} \sum_{j=N_{K-1}+1}^{N_K}|E_1(j)|^p \Big]^{1/p} &\le \kappa_1 \gamma_{n_{N_{{K-1}}}+1}\beta^{n_{N_{{K-1}}}+1}=: \mathcal E_1(K) \\
		\E \Big[\1_{\{|\theta_{n_{N_{K-1}}}| \le R\}} \sum_{j=N_{K-1}+1}^{N_K}|E_i(j)|^p \Big]^{1/p} &\le\kappa_i \mathcal E_i(K) \quad i=2,3,4, \\
	\end{align*}
	where $\mathcal E_i(K)$ for $i=2,3,4$ is defined in \eqref{eq:mathcalE} with $\ell(K)=N_{K-1}$ and $r(K)=N_K$.
	In order to bound the last term in \eqref{eq:2930472042}, note that on $\{|\theta_{n_{N_{K-1}}}|\le R\}$ and for all $j \in \{N_{K-1}+1, \dots, N_K\}$ and $0 \le t \le t_{n_{j}}-t_{n_{j-1}}$ one has almost surely that $|\Psi_{t}(\theta_{n_{j-1}})-\theta_{n_{j-1}}|\le C_1 t$ and, thus, 
	\begin{align*}
		|(t_{n_j}-t_{n_{j-1}})f(\theta_{n_{j-1}})-\Psi_{t_{n_j}-t_{n_{j-1}}}(\theta_{n_{j-1}})| &\le \int_{0}^{t_{n_j}-t_{n_j-1}} |f(\theta_{n_{j-1}})- f(\Psi_t(\theta_{n_{j-1}})) | \, dt   \\
		& \le \int_{0}^{t_{n_j}-t_{n_j-1}} C_2 C_1 t \, dt = \frac 12 C_2 C_1 (t_{n_j}-t_{n_{j-1}})^2.
	\end{align*}
	Therefore, there exists a $\kappa_5>0$ such that
	\begin{align*}
		\E \Big[\1_{\{|\theta_{n_{N_{K-1}}}| \le R\}} \sum_{j=N_{K-1}+1}^{N_K}|E_5(j)|^p \Big]^{1/p} &\le \kappa_5 \mathcal E_2(j)
	\end{align*}
	and we define $\mathcal E_5(j):=\mathcal E_2(j)$.
	Now, using Markov's inequality and the triangle inequality for the $L^p$-norm 
	\begin{align}\begin{split} \label{eq:23965234}
			&\P \Bigl(\Bigl\{ \1_{\{|\theta_{n_{N_{K-1}}}|<R\}} \sup_{\ell=N_{K-1}, \dots, N_{K}} \sup_{i=\ell+1, \dots, N_{K}}  |\theta_{n_i} - \Psi_{t_{n_i}-t_{n_\ell}}(\theta_{n_\ell})| > \eps\Bigr\} \Bigr)^{1/p} \\
			&\le \frac{1}{\eps^p} \E\Bigl[\1_{\{|\theta_{n_{N_{K-1}}}|<R\}}\Bigl( \sup_{\ell=N_{K-1}, \dots, N_{K}} \sup_{i=\ell+1, \dots, N_{K}} |\theta_{n_i} - \Psi_{t_{n_i}-t_{n_\ell}}(\theta_{n_\ell})|\Bigr)^p \Bigr]^{1/p} \\
			& \le  \frac{1}{\eps^p} \sum_{i=1}^5 \kappa_i \mathcal E_i(K).
		\end{split}
	\end{align}
	Since Proposition~\ref{prop:summable} implies that $\sum_{K \in \N}\mathcal E_i^p(K)<\infty$ for all $i=1, \dots, 5$, we get by Borel-Cantelli's lemma that almost surely
	$$
	\1_{\{|\theta_{n_{N_{K-1}}}|<R\}} \sup_{i=\ell+1, \dots, N_{K}} \sup_{\ell=N_{K-1}, \dots, N_{K}} |\theta_{n_i} - \Psi_{t_{n_i}-t_{n_\ell}}(\theta_{n_\ell})| < \eps
	$$
	for all but finitely many $K$'s. This proves \eqref{eq:2834656252}.
	
	To statement \ref{thm:item_ii}, note that we can still use \eqref{eq:23965234} so that the statement follows if one can show that $\mathcal E_i(K) \to 0$ for all $i=1, \dots, 5$.
	$\mathcal E_1(K)\to 0$ is due to $\gamma_n \to 0$. Using \eqref{eq:238969000}, we have
	$$
	\mathcal E_2(K) \le \sum_{i=n_{N_{K-1}}+1}^{n_{N_K}} \rho_0 \gamma_k^{4/3} \le T \gamma_{n_{N_{K-1}}+1}^{1/3} \to 0
	$$
	and, thus, $\mathcal E_3(K)\to 0$ and $\mathcal E_5(K)\to 0$. Similarly, using \eqref{eq:23975843},
	$$
	\mathcal E_4(K) \lesssim \gamma_{n_{N_{K-1}}+1}^{1/3} \to 0.
	$$

\section{Proof of Theorem~\ref{thm:degeneratenoiseNEW}} \label{sec:proof2}

First, we show a general capture statement for the semiflow~\eqref{eq:ODE}.
\begin{lemma} \label{lem:capture}
	Consider the setup of \cref{thm:item_ii} in \cref{thm:main}.
	Let $V \subset \mathcal D \subset \R^d$ be open and $K \subset \R^d$ be a compact set, and assume that there exist $T, \eps >0$ such that for all $\theta \in V\cap K$ with $\sup_{t\in[0,T]} d(\Psi_t(\theta),K)\le \eps $ it holds that
	\bas{\label{eq:6344}
		\bigcup_{0\le t \le T} B_\eps(\Psi_t(\theta)) \subset \mathcal D \quad \text{ and } \quad B_\eps(\Psi_T(\theta)) \subset V.
	}
	Then, for almost all 
	$$
	\omega \in  \{\theta_n \in V \text{ for infinitely many } n \in \N_0\} \cap  \bigcup_{n_0\in\N_0}\bigcap_{n \ge n_0}\{\theta_n \in K\}
	$$
	one has $\theta_n(\omega) \in \mathcal D$ for all but finitely many $n$.
\end{lemma}

\begin{proof}
	We fix an $n_0 \in \N_0$ and a
	$$
	\omega \in  \{\theta_n \in V \text{ for infinitely many } n \in \N_0\} \cap \bigcap_{n \ge n_0}\{\theta_n \in K\}
	$$
	such that the conclusion of part  \ref{thm:item_ii}  of \cref{thm:main} holds that is
	\begin{align*}
		\lim_{N\to\infty} \sup_{\substack{n\in\N_0: \\ 0 \le t_n-t_N \leq T}}\left| \theta_{n}(\omega)-\Psi_{t_n-t_N}( \theta_{N}(\omega))\right|  =0.
	\end{align*}
	Recall that the \Adam\ vector field is uniformly bounded, see Lemma~\ref{lem:lipschitz}. 
	We let $N_0 \ge n_0$ such $\|f\|_\infty \gamma_{N_0}\le \eps/2$ and 
	for all $N\ge N_0$,
	$$
	\sup_{\substack{n\in\N_0: \\ 0 \le t_n-t_N \leq T}}\left| \theta_{n}(\omega)-\Psi_{t_n-t_N}( \theta_{N}(\omega))\right|  \le \frac \eps2.
	$$
	This implies that $\sup_{t\in[0,T]} d(\Psi_t(\theta_N(\omega)),K)\le \eps$.
	We show that, if $\theta_{N_1}(\omega)$ is in $V$ for an $N_1\ge N_0$, it will never leave $\cD$ again. Indeed, one has for $n=N_1,\dots,N_2:=\max\{n\in\N: t_{n}-t_{N_1}\le T\}$ that
	$$
	|\theta_n(\omega)-\Psi_{t_n-t_{N_1}}(\theta_{N_1}(\omega))|\le \frac\eps2
	$$ 
	which entails with the left-hand side of (\ref{eq:6344}) and the fact that $\theta_{N_1}(\omega)\in V\cap K$ that $(\theta_n(\omega))_{n=N_1,\dots,N_2}$ stays in $\cD$. Moreover, 
	\begin{align*}
		|&\theta_{N_2}(\omega)-\Psi_{T}(\theta_{N_1}(\omega))|\\
		&\le 
		|\theta_{N_2}(\omega)-\Psi_{t_{N_2}-t_{N_1}}(\theta_{N_1}(\omega))|
		+ |\Psi_{t_{N_2}-t_{N_1}}(\theta_{N_1}(\omega))-\Psi_{T}(\theta_{N_1}(\omega))|\\
		&\le 	\frac\eps2 +\frac\eps2=\eps
	\end{align*}
	so that since $\theta_{N_1}(\omega)\in V\cap K$ and $B_\eps(\Psi_T(\theta_{N_1}(\omega))\cap K \subset V$ we get that $\theta_{N_2}(\omega)\in V\cap K$. Now one iterates the argument and obtains that $(\theta_n(\omega))_{n\ge N_1}$ will stay in $\cD$.
\end{proof}

Next, we show that, close to the set of critical points, $R$ is a Lyapunov function for the semiflow \eqref{eq:ODE}. For this, let us cite Lemma~8.5 in \cite{dereich2024convergence}.

\begin{lemma} \label{lem:dereich2}
	Let $\beta, \delta, q \in (0,\infty)$ with $q<\sqrt{\beta}<1$ and let $\mathcal Z:=(Z_k)_{k \in -\N_0}$ be an i.i.d. sequence of real-valued random variables satisfying for all $k \in -\N_0$ that
	\begin{align*}
		\P(Z_k^2 < \delta) \le q.
	\end{align*}
	Then,
	\begin{align*}
		\E \Bigl[ \frac{1}{\sqrt{v(\mathcal Z)}}\Bigr] \le  \Bigl(\frac{\beta}{1-\beta}\Bigr)^{\frac 12}  \frac{1-q}{\beta^{\frac 12 }-q}\delta^{-\frac 12 },
	\end{align*}
	where
	$$
	v: \R^{-\N_0}  \to [0,\infty] \quad, \quad (x_k)_{k \in -\N_0} \mapsto (1-\beta) \sum_{k \in -\N_0} \beta^{-k} x_k^2.
	$$
\end{lemma}

\begin{lemma} \label{lem:scalar}
	Let $R: \R^d \to [0,\infty)$ be a $C^1$-function such that $\E[X(U, \theta)]=-\nabla R(\theta)$ for all $\theta \in \R^d$. 
	Let $\mathcal D\subset \R^d$ and assume that there exist constants $C_1, C_2, p_1,p_2,q\in (0,\infty)$ with $1< p_2\le 2$,  $\sfrac{p_2}{2}+1<p_1\le \frac 32 p_2$ and $q<\sqrt{\beta}$ such that for all $\theta \in \mathcal D$ and $i=1, \dots, d$
	\begin{align} \label{eq:9982622}
		\E[|X^{(i)}(U, \theta)|^3] \le C_1 |\nabla R(\theta)|^{p_1} \quad \text{ and } \quad \P(|X^{(i)}(U,\theta)|^2< C_2 |\nabla R^{(i)}(\theta)|^{p_2})\le   q.
	\end{align}
	Then there exists a $\delta\in(0,\infty)$, that only depends on $C_1, C_2, p_1,p_2,q, \alpha, \beta, \epsilon$, such that one has for all $\theta \in \mathcal D \cap \{\theta' \in \R^d: |\nabla R(\theta')|\le \delta \}$  and the respective \Adam\ vector field $f$ that
	$$
	\langle f(\theta), \nabla R(\theta) \rangle < 0 \quad \text{ or } \quad  f(\theta)= \nabla R(\theta)=0.
	$$
\end{lemma}

\begin{proof}
	We recall the definition of the \Adam\ vector field for $i=1, \dots, d$
	$$
	f^{(i)}(\theta) = (1-\alpha) \E\Biggl[\Biggl( \sqrt{(1-\beta)\sum_{k=-\infty}^0 \beta^{-k} X^{(i)}(U_k,\theta)^2}+\epsilon \Biggr)^{-1} \sum_{k=-\infty}^0 \alpha^{-k} X^{(i)}(U_k,\theta) \Biggr].
	$$
	Due to the first assumption in \eqref{eq:9982622}, for $\theta \in \mathcal D$ with $\nabla R(\theta) =0$ one has $X^{(i)}(U,\theta)=0$ almost surely for all $i \in \{1, \dots, d\}$ so that $f(\theta)=0$. 
	
	Now, let $\theta \in \mathcal D$ and assume that $\nabla R(\theta)\neq 0$. 
	Then,
	\begin{align*}
		\langle f(\theta), \nabla R(\theta) \rangle =\sum_{i=1}^d \1_{\{\nabla R^{(i)}(\theta) \neq 0\}} f^{(i)}(\theta) \, \nabla R^{(i)}(\theta).
	\end{align*}
	Let $i\in \{1, \dots, d\}$ with $\nabla R^{(i)}(\theta) \neq 0$.
	We will approximate $f^{(i)}$ by the function $\tilde f^{(i)}: \R^d \to \R$ given by 
	$$
	\tilde f^{(i)}(\theta) = (1-\alpha) \sum_{k=-\infty}^0 \alpha^{-k} \E\Biggl[ \Bigl(\sqrt{ G_{\neq k}^{(i)}(\theta)}+\epsilon\Bigr)^{-1} \Biggr] \E[X^{(i)}(U_k, \theta)],
	$$	
	where
	$$
	G_{\neq k}^{(i)}(\theta) := (1-\beta) \sum_{r \in -\N_0: r \neq k } \beta^{-r} X^{(i)}(U_r, \theta)^2.
	$$
	For this approximation we get
	\begin{align*}
		\tilde f^{(i)}(\theta)\,  \nabla R^{(i)}(\theta)  \le -  \E\Biggl[ \Biggl(\sqrt{ (1-\beta) \sum_{r \in -\N_0 } \beta^{-r} X^{(i)}(U_r, \theta)^2}+\epsilon\Biggr)^{-1} \Biggr] |\nabla R^{(i)}(\theta)|^2.
	\end{align*}
	Using Jensen's inequality and the first assumption in \eqref{eq:998262} we get
	\begin{align*}
		\E\Biggl[ \Biggl(\sqrt{ (1-\beta) \sum_{r \in -\N_0 } \beta^{-r} X^{(i)}(U_r, \theta)^2}+\epsilon\Biggr)^{-1} \Biggr] &\ge \Biggl(\E\Biggl[ \sqrt{ (1-\beta) \sum_{r \in -\N_0 } \beta^{-r} X^{(i)}(U_r, \theta)^2}\Biggr] +\epsilon  \Biggr)^{-1} \\
		&\ge \Biggl(\sqrt{\E\Biggl[  (1-\beta) \sum_{r \in -\N_0 } \beta^{-r} X^{(i)}(U_r, \theta)^2\Biggr]} +\epsilon  \Biggr)^{-1} \\
		& \ge \Bigl(\E[  |X^{(i)}(U_r, \theta)|^3]^{\frac 13} +\epsilon  \Bigr)^{-1} \ge  \Bigl(C_1^{\frac 13} |\nabla R(\theta)|^{\frac{p_1}{3}} +\epsilon  \Bigr)^{-1}
	\end{align*}
	so that 
	\begin{align} \label{eq:9068033}
		 \tilde f^{(i)}(\theta) \,  \nabla R^{(i)}(\theta) \le - |\nabla R^{(i)}(\theta)|^2  \Bigl(C_1^{\frac 13} |\nabla R(\theta)|^{\frac{p_1}{3}} +\epsilon  \Bigr)^{-1}.
	\end{align}
	
	Next we compute $|f^{(i)}(\theta)-\tilde f^{(i)}(\theta)|$.
	Let $h(x)=\frac{1}{\sqrt{x}+\eps}$ and note that
	$$
	-h'(x)= \frac{1}{(\sqrt x +\eps)^2} \frac{1}{2\sqrt x}.
	$$
	is monotonically decreasing. Using the mean value theorem, for $0<x\le y$ there exists a $\xi \in [x,y]$ such that
	$$
	h(x)-h(y) = h'(\xi) (x-y) \le h'(x) (y-x).
	$$
	Setting  $x= G_{\neq k}^{(i)}(\theta)$ and $y = (1-\beta) \sum_{r=-\infty}^0 \beta^{-r}X^{(i)}(U_r, \theta)^2$ we get that
	\bas{\label{eq:23476}
		&|f^{(i)}(\theta)-\tilde f^{(i)}(\theta)| \le -(1-\alpha) \sum_{k=-\infty}^0 \E \Biggl[ \alpha^{-k} \, |X^{(i)}(U_k, \theta)| \,  h'(G_{\neq k}^{(i)}(\theta)) \,  (1-\beta) \beta^{-k} X^{(i)}(U_k, \theta)^2  \Biggr] \\
		& = (1-\alpha)(1-\beta) \sum_{k=-\infty}^0 (\alpha \beta)^{-k} \,  \E \left[  \frac{1}{2\Bigl(\sqrt{G_{\neq k}^{(i)}(\theta)}+\epsilon\Bigr)^2\sqrt{G_{\neq k}^{(i)}(\theta)}} \right] \E[|X^{(i)}(U_k, \theta)|^3].
	}
	
	Since $\nabla R(\theta)^{(i)} \neq 0$, we can use the second inequality in \eqref{eq:9982622} together with Lemma~\ref{lem:dereich2} to deduce that there exists a constant $\kappa_1 \in(0,\infty)$ only depending on $q,\beta$ and $C_2$  such that 
	\begin{align*}
		\E \left[ \Biggl( 2\Bigl(\sqrt{G_{\neq k}^{(i)}(\theta)}+\epsilon\Bigr)^2\sqrt{G_{\neq k}^{(i)}(\theta)}\Biggr)^{-1} \right] \le \E \left[ \Bigl( 2\eps^2\sqrt{G_{\neq k}^{(i)}(\theta)}\Bigr)^{-1} \right] 
		&\le \E \left[ \Bigl( 2\eps^2\sqrt{G_{\neq 0}^{(i)}(\theta)}\Bigr)^{-1} \right]  \\
		& \le  \frac{\kappa_1}{2\eps^2} |\nabla R(\theta)^{(i)}|^{-\sfrac{p_2}{2}}.
	\end{align*}
	
	Together with (\ref{eq:23476}) and  the first inequality in \eqref{eq:9982622} we conclude that there exists a constant $\kappa_2 \in( 0,\infty)$ only depending on $q,\alpha, \beta,\delta,\eps, C_1$ and $C_2$ such that	
	\begin{align} \label{eq:1}
		|f^{(i)}(\theta)-\tilde f^{(i)}(\theta)| \le \kappa_2 |\nabla R(\theta)|^{p_1} \,  |\nabla R(\theta)^{(i)}|^{-\sfrac{p_2}{2}}.
	\end{align}	
	Thus, since  $\frac{p_2}{2}\le 1$ we can combine \eqref{eq:9068033} and \eqref{eq:1} to get
	\begin{align} \begin{split}
			\label{eq:2}
			&\langle f(\theta) ,  \nabla R(\theta)\rangle \\
			  &\le  \sum_{i=1}^d \1_{\{\nabla R^{(i)}(\theta) \neq 0\}} \tilde f^{(i)}(\theta) \, \nabla R^{(i)}(\theta) + \sum_{i=1}^d \1_{\{\nabla R^{(i)}(\theta) \neq 0\}} |f^{(i)}(\theta)-\tilde f^{(i)}(\theta)| \, |\nabla R^{(i)}(\theta)| \\
			& \le - |\nabla R(\theta)|^2  \Bigl(C_1^{\frac 13} |\nabla R(\theta)|^{\frac{p_1}{3}} +\epsilon  \Bigr)^{-1} + d \kappa_2 |\nabla R(\theta)|^{p_1} \max_{i=1, \dots, d}  \1_{\{\nabla R^{(i)}(\theta) \neq 0\}} |\nabla R^{(i)}(\theta)|^{1-\frac{p_2}{2}}\\
			& \le - |\nabla R(\theta)|^2  \Bigl(C_1^{\frac 13} |\nabla R(\theta)|^{\frac{p_1}{3}} +\epsilon  \Bigr)^{-1} +  d \kappa_2 |\nabla R(\theta)|^{p_1+1- \frac{p_2}{2}}.
		\end{split}
	\end{align}	
	Since $1+p_1-\sfrac{p_2}{2} > 2$ there exists a $\delta>0$, that only depends on $C_1, C_2, p_1,p_2,q, \alpha, \beta, \epsilon$, such that, for all $\theta \in V \cap \{\theta' \in \R^d: |\nabla R(\theta')|\le \delta \}$,
	$$
	\langle f(\theta), \nabla R(\theta) \rangle  < 0 .
	$$
\end{proof}

\begin{proof} [Proof of Theorem~\ref{thm:degeneratenoiseNEW}]
	Without loss of generality, we assume that $\arg\min R \big|_{\mathcal D}=0$.
	For $\eps>0$, we write $K_\eps:=\{ \theta \in \R^d: d(\theta, K) \le \eps\}$. 
	
	By definition, the open set $\cD$ satisfies $\mathcal D \supset \mathcal C \cap K$ and first we show that there exists $\eps\in(0,\infty)$ such that also $\mathcal D \supset \mathcal C \cap K_\eps$. If this were not true, then there exist a sequence $(\eps_n)_{n \in \N}$ with $\eps_n \to 0$ and a convergent sequence $(x_n)_{n \in \N}$ with 
	$x_n \in (\mathcal C \cap K_{\eps_n})\backslash \cD$ for all $n\in\N$. By compactness, this sequence would possess an accumulation point $x\in \cC\cap K$. Since $\cD^c$ is closed this point would not lie in $\cD$ contradicting $\cC\cap K\subset \cD$.
	
	Using Lemma~\ref{lem:scalar} together with \eqref{eq:998262}, there exist $\delta' >0$ such that for all $\theta \in K_\eps \cap \mathcal D$ with $|\nabla R(\theta)|<\delta'$ one has
	\begin{align} \label{eq:763546}
		\langle f(\theta), \nabla R(\theta) \rangle < 0 \quad \text{ or } \quad  f(\theta)= \nabla R(\theta)=0.
	\end{align}
	
	Recall that by assumption (\ref{eq:783256}) one has for every $\theta\in\overline\cD$ that $\nabla R(\theta)=0$ iff $\theta\in\cC$ and also $R(\theta)=0$ iff $\theta\in\cC$. 
	After possibly shrinking $\mathcal D$ we get to an open neighbourhood of $\cC\cap K_\eps$ satisfying$\sup_{\theta \in \mathcal D} |\nabla R(\theta)|<\delta'$.
	Now there exists a $\delta \in (0, \infty)$ such that 
	\bas{
		\{\theta \in K_\eps \cap \overline{\mathcal D}: R(\theta) \le \delta\}  \subset \mathcal D.
	}
	Otherwise, there existed a $K_\eps \cap \partial \cD$-valued sequence $(x_n)_{n \in \N}$ with  $\lim_{n\to\infty} R(x_n) =0$. Since $R$ is continuous the limit point $x$ would satisfy $x\in K_\eps\cap \partial \cD$ and $R(x)=0$ and, hence, $x\in\cC\cap K_\eps\subset \cD$  with $\cD\cap \partial D=\emptyset$ which is a contradiction.

	Next, we want to apply Lemma \ref{lem:capture} with the sets 
	\bas{
		V:=\{\theta \in \R^d: R(\theta)<\delta/2\} \text{ \ \ and \ \ }  \bar V := \{\theta \in \R^d: R(\theta)<\delta\}.}
	For this, consider the compact set $W:=\{\theta \in K_{\eps}:\delta/4\le  R(\theta)\le \delta\}$. By property (\ref{eq:763546}) we have that
	\bas{
		-\eta:=\max_{\theta\in W} \langle f(\theta),\nabla R(\theta)\rangle<0.}
	For a fixed $\theta \in K\cap V$ we let $\tau:=\inf\{s \ge 0 : \Psi_s(\theta) \notin  K_{\eps}\cap \mathcal D \vee R(\Psi_s(\theta)) > \delta\}$ and observe that one has $\theta \in \mathcal C$ or for all $t\in[0,\tau]$ that
	$$
	R(\Psi_t(\theta)) = R(\theta) + \int_0^t \langle f(\Psi_s(\theta)), \nabla R(\Psi_s(\theta)) \rangle \, ds <  R(\theta) \le \delta.
	$$
	In the case where $\tau$ is finite we have that $R(\Psi_\tau(\theta))<\delta$ and, thus, $\Psi_\tau(\theta) \in \mathcal D$. This entails that $d(\Psi_\tau(\theta), K) \ge \eps$.
	
	Note that  $R$ is uniformly continuous on the compact set $K_{\eps}$ and we can pick $\eps'\in(0,\eps)$  so that
	\begin{align}\label{eq83546}
		\sup_{\theta,\theta'\in K_{\eps'}:|\theta-\theta'|\le \eps'}  |R(\theta') - R(\theta)| < \frac{\delta}{8}.
	\end{align}
	Choose $T= \frac {\delta}{4\eta}$ and observe that  for arbitrarily fixed  $\theta \in V \cap K$ with $\sup_{t \in [0,T]} d(\Psi_t(\theta),K)\le \eps'$ one has that the path $(\Psi_t(\theta))_{t\in[0,T]}$ stays in $K_{\eps}$ and satisfies  
	\begin{align*}
		R(\Psi_t(\theta)) &= R(\theta) + \int_0^t \langle f(\Psi_s(\theta)), \nabla R(\Psi_s(\theta)) \rangle \, ds \\
		&\le R(\theta) - \int_0^t \1_{\{R(\Psi_s(\theta))\in [\delta/4,\delta]\} }\eta \, ds\le \max(R(\theta)-\eta t, \delta/4).
	\end{align*}
	In particular, $R(\Psi_t(\theta))\le \delta/2$ for all $t\in[0,T]$ and $R(\Psi_T(\theta))\le \delta/4$ so that (\ref{eq83546}) entails that
	\bas{
		\bigcup_{0 \le t \le T} B_{\tilde \eps}(\Psi_t(\theta))\subset  \bar V \quad \text{ and } \quad B_{\tilde \eps}(\Psi_T(\theta)) \subset V.
	}
	Using Lemma~\ref{lem:capture}, we get for almost all
	$$
	\omega \in \bigl\{\theta_n \in V \text{ for infinitely many } n \in \N_0\bigr\} \cap  \bigcup_{n_0\in\N_0}\bigcap_{n \ge n_0}\{\theta_n \in K\} =: \Omega_0
	$$
	that $\theta_n(\omega) \in \overline V$ for all but finitely many $n$. 
	Now, $\mathcal L((\theta_n(\omega))_{n \in \N_0})$ is closed and invariant for the semiflow \eqref{eq:ODE}, see e.g. \cite[Theorem~0.1]{benaim1996asymptotic}, and contained in $\overline V \cap K$ such that $R$ is a Lyapunov function for $\mathcal C \cap \mathcal L((\theta_n(\omega))_{n \in \N_0})$. Thus, we can apply Corollary~\ref{cor:Lyapunov}, to deduce that $d(\theta_n,\mathcal C)\to 0$ almost surely on $\Omega_0$.

\end{proof}

\section{Proof of Theorem~\ref{thm:ERM}} \label{sec:proof3}
As a key auxiliary step in the proof Theorem~\ref{thm:ERM}, we establish the following combinatorial result.

\begin{theorem} \label{thm:combinatoric}
	Let $N,M \in \N$ and denote $\mathcal X^N=\{x \in \R^N : \sum_{i=1}^N x_i = 1\}$. Let $I_1, \dots, I_M$ be independent, $\mathrm{Unif}(\{1, \dots, N\})$-distributed random variables. Then, there exists a $q_0>0$ such that for all $0\le q < q_0$ one has that
	$$
	\sup_{x\in \mathcal X^N}\P\Big( \frac{1}{M} \Big| \sum_{k=1}^M x_{I_k} \Big|\le q  \Big) \le \frac 12 \Big(1-\frac{1}{N}\Big)^M+ \frac 12.
	$$
\end{theorem}

The proof of Theorem~\ref{thm:combinatoric} is divided into two lemmas.

\begin{lemma}
	Let $N,M \in \N$ and denote $\mathcal X^N=\{x \in \R^N : \sum_{i=1}^N x_i = 1\}$ and $\mathcal Z^N= \{z \in \R^N : \|z\|_{\infty} =1\}$. Let $I_1, \dots, I_M$ be independent, $\mathrm{Unif}(\{1, \dots, N\})$-distributed random variables. Then, there exists a $q_0>0$ such that for all $0\le q < q_0$ one has
	$$
	\sup_{x \in \mathcal X^N}\P\Big( \frac{1}{M} \Big| \sum_{k=1}^M x_{I_k} \Big|\le q  \Big) \le \sup_{z \in \mathcal Z^N}\P\Big( \frac{1}{M} \sum_{k=1}^M z_{I_k}  = 0 \Big).
	$$
\end{lemma}

\begin{proof}
	Consider the random vector $C=(C_1, \dots, C_n)^\dagger$ with  
	$$
	C_i = \sum_{j=1}^M \1_{\{I_j=i\}} \text{, \ \ for every } i=1,\dots,N.
	$$
	$C$ takes values in the finite set 
	\begin{align*}
		\mathcal C = \Bigl\{c \in \N_0^N: \sum_{i=1}^N c_i=M\Bigr\}
	\end{align*}	
	and for every $c \in \cC$  we denote by $p_c = \P(C=c)$ the respective probability. Note that, for every $x \in \mathcal X^N$, one has
	$$
	\P\Big( \frac{1}{M} \Big| \sum_{k=1}^M x_{I_k} \Big|\le q  \Big) =\sum_{c \in \mathcal C} p_c \1_{\{|\langle c,x\rangle | \le Mq\}}.
	$$
	
	We will relate the latter probability to the occurrence of so called zero patterns:  for every $z \in \mathcal Z^N$ we define its zero-pattern by
	$$
	\mathcal C_z:= \{c \in \mathcal C: \langle c, z \rangle =0\} \subsetneq \mathcal C.
	$$
	Recall that $\cC$ is finite and and note that the mapping 
	$\R^N \ni z \mapsto (\langle c,z \rangle)_{c\in\cC}\in\R^{\cC}$ is continuous.
	Thus, for every $z \in \mathcal Z^N$, there exists a neighborhood $U_z \subset \R^N$ and a constant $\delta_z >0$ such that for every $z' \in U_z$ and every $c \notin \mathcal C_z$ one has $|\langle c,z' \rangle |\ge \delta_z$.
	By compactness of $\mathcal Z^N$, we can choose a finite set $\{z^{(1)}, \dots, z^{(k)}\}$ with $\bigcup_{i=1, \dots, k} U_{z^{(i)}} \supset \mathcal Z^N$ and we set $\delta^* = \min_{i=1, \dots, k} \delta_{z^{(i)}}$.
	
	Note that for every $x \in \mathcal X^N$ one has $\|x\|_\infty \ge 1/N$ and $x/\|x\|_\infty \in \mathcal Z^N$. Thus, there exists an $i \in \{1, \dots, k\}$ with $x/\|x\|_\infty \in U_{z^{(i)}}$. This entails that  for all $c \notin \mathcal C_{z^{(i)}}$ one has that
	$$
	|\langle c,x \rangle | \ge \|x\|_\infty \,  \Big|\Big\langle c, \frac{x}{\|x\|_\infty} \Big\rangle \Big| \ge \frac {\delta^*}N .
	$$
	Set $q_0 = \delta^*/(NM)$ and note that for all $q<q_0$
	$$
	\P\Big( \frac{1}{M} \Big| \sum_{k=1}^M x_{I_k} \Big|\le q  \Big) =\sum_{c \in \mathcal C} p_c \,\1_{\{|\langle c,x\rangle | \le Mq\}} \le \sum_{c \in \mathcal C} p_c \,\1_{\{|\langle c,z^{(i)}\rangle | = 0 \}} = \P\Big( \frac{1}{M} \sum_{k=1}^M z_{I_k}  = 0 \Big).
	$$
	Taking the supremum over all $x \in \mathcal X^N$ finishes the proof. 
\end{proof}

\begin{lemma}
	Let $N,M \in \N$, $z=(z_1, \dots, z_N) \in \R^N$ with $z \neq 0$. Let $I_1, \dots, I_M$ be independent, $\mathrm{Unif}(\{1, \dots, N\})$-distributed random variables. Then,
	$$
	\P\Big(\sum_{k=1}^M z_{I_k} = 0 \Big) \le \frac{1}{2}\Big(1-\frac 1N\Big)^M + \frac 12.
	$$
\end{lemma}

\begin{proof}
	Set $p=\P(z_{I_1}\neq 0)$ and $\tilde z =(\tilde z_1, \dots, \tilde z_\tau)$ be a tuple that contains the non-zero entries in $z$. Then,
	\begin{align*}
		\P\Big(\sum_{k=1}^M z_{I_k}  = 0\Big) = (1-p)^M + \sum_{k=1}^M {M \choose k} p^k (1-p)^{M-k} \alpha_k,
	\end{align*}
	where for $k \in \N$ and $J_1, \dots, J_k$ independent, $\mathrm{Unif}(\{1, \dots, \tau\})$-distributed random variables we define
	$$
	\alpha_k = \P\Big(\sum_{i=1}^k \tilde z_{J_i}=0\Big).
	$$
	We will show that  $\alpha_k \le \frac 12$ for all $k \in \N$. Since $p\ge \frac 1N$ this implies
	$$
	\P\Big(\sum_{k=1}^M z_{I_k}  = 0\Big) \le (1-p)^M + \frac{1}{2} (1-(1-p)^M) \le \frac{1}{2}\Big(1-\frac 1N\Big)^M + \frac 12,
	$$
	where we have used that the second term is monotonically decreasing in $p$.

	We are left with proving the following claim:
	Let $\tau, k \in \N$ and $\tilde z_1, \dots, \tilde z_\tau \neq 0$. Then,
	$$
	\P\Big(\sum_{i=1}^k \tilde z_{J_i}=0\Big) \le \frac 12,
	$$
	where $J_1, \dots, J_k$ are independent, $\mathrm{Unif}(\{1, \dots, \tau\})$-distributed random variables.
	
	We split the indices into $P := \{i: \tilde x_i > 0\}$ and $N := \{ i: \tilde x_i < 0  \}$. Without loss of generality, one has $|N|\ge |P|$. Thus, we can define an injective mapping $\varphi:P \to N$. We define an equivalence relation on the set $\{1, \dots, \tau\}^k$. Let $(i_1, \dots, i_\tau) \in \{1, \dots, \tau\}^k$. Denote by $j$ the smallest index in $\{1, \dots, \tau\}$ such that $i_j \in P \cup \varphi(N)$ and $j=\infty$ if there is no such index. We say that $(i_1, \dots, i_\tau) \sim (\tilde i_1, \dots, \tilde i_\tau) $ if and only if 
	$$
	\begin{cases}
		\tilde i_k = i_k, & \text{ if } k \neq j,\\
		\tilde i_j = \varphi(i_j), & \text{ if } i_j \in P, \\
		\tilde i_j = \varphi^{-1}(i_j), & \text{ if } i_j \in \varphi(P),
	\end{cases}
	$$
	and $(i_1, \dots, i_\tau) = (\tilde i_1, \dots, \tilde i_\tau)$ if $j = \infty$. Note that every equivalence class has at most two elements and at least one of the elements leads to a sum which is non-negative. We have, thus, proved the claim.
\end{proof}

\begin{proof}[Proof of Theorem~\ref{thm:ERM}]
	The proof follows from Theorem~\ref{thm:degeneratenoiseNEW} if we can verify the assumptions in \eqref{eq:998262} for the innovation $(X_M, \mathbf U)$ defined in \eqref{eq:ERMmini-batch} and the choice $p_1=3, p_2=2$.
	Regarding the first assumption in \eqref{eq:998262}, note that
	\begin{align*}
		\E[&|\nabla R_U(\theta)|^3 ] =  \frac 1N \sum_{i=1}^N  |(\mathfrak N^{\theta}(y_i)-z_i)\grad_\theta \mathfrak N^{\theta}(y_i)|^3.
	\end{align*}
	After possibly shrinking $U$, we can assume that $U$ is bounded so that the gradient $\grad_\theta \mathfrak N^\theta (y_i)$ is bounded on $U$ for all $i=1, \dots, N$. Thus, there exists a constant $\kappa(N) \ge 0$ such that 
	$$
	\E[|\nabla R_U(\theta)|^3 ] \le \kappa(N) R(\theta)^{\frac 32} \quad \text{ for all }\theta \in U.
	$$
	Since $R$ satisfies the PL-inequality on $U$, i.e. there exists a $\mu \ge 0$ such that for all $\theta \in U$
	$$
	2 \mu R(\theta) \le |\nabla R(\theta)|^{2},
	$$
	one has
	\begin{align}\label{eq:moment}
		\E[|\nabla R_U(\theta)|^3 ] \le \kappa(N) R(\theta)^{\frac 32}\le \kappa(N) (2\mu) ^{-\frac 32} |\nabla R(\theta)|^3,
	\end{align}
	which verifies the first assumption in \eqref{eq:998262}.
	Regarding the second assumption, fix $\theta \in U$ and $i \in \{1, \dots, d\}$ and assume that $\nabla R^{(i)}(\theta) \neq 0$. Then,
	$x=(x_1, \dots, x_N) \in \R^d$ given by $x_j = \nabla R_j^{(i)}(\theta)/(N \, |\nabla R^{(i)}(\theta)|)$ satisfies 
	$$
		\sum_{j=1}^N x_j = \frac{\nabla R^{(i)}(\theta)}{|\nabla R^{(i)}(\theta)|} \in \{\pm 1\}.
	$$
	Thus, Theorem~\ref{thm:combinatoric} shows that there exists a $q>0$ such that
	\begin{align*}
		\P(|X_M^{(i)}(\mathbf U,\theta)|^2< q |\nabla R^{(i)}(\theta)|^{2}) & = \P\Big(\frac{1}{M}\Big|\sum_{j=1}^M\nabla R_{U_j}^{(i)}(\theta)\Big|^2< q |\nabla R^{(i)}(\theta)|^{2}\Big)\\
		& = \P\Big( \frac{1}{M} \Big| \sum_{k=1}^M x_{I_k} \Big|\le Nq  \Big) \le \frac 12 \Big(1-\frac{1}{N}\Big)^M+ \frac 12\le \sqrt \beta,
	\end{align*}
	where in the last step we used the assumption on $M$. This verifies the second assumption in \eqref{eq:998262}.
\end{proof}

\subsection*{Acknowledgments}
SK acknowledges funding by the Deutsche Forschungsgemeinschaft (DFG, German Research Foundation) – CRC/TRR 388 ''Rough Analysis, Stochastic Dynamics and Related Fields'' – Project ID 516748464.
Moreover, this work has been partially funded by the European Union (ERC, MONTECARLO, 101045811). 
The views and the opinions expressed in this work are however those of the authors only and do not 
necessarily reflect those of the European Union or the European Research Council (ERC). 
Neither the European Union nor the granting authority can be held responsible for them. 
Furthermore, this work has been funded by the Deutsche Forschungsgemeinschaft (DFG, German Research Foundation) under Germany's Excellence Strategy EXC 2044-390685587, Mathematics Münster: Dynamics-Geometry-Structure.

\bibliographystyle{alpha}
\bibliography{Adamnew}

\end{document}